\documentclass[final,3p,times]{elsarticle}

\usepackage{graphicx}
\usepackage[ansinew]{inputenc}
\usepackage{amssymb}
\usepackage{latexsym}
\usepackage{amsfonts}
\usepackage{amsmath}
\usepackage{color}
\usepackage{subfigure}
\usepackage{lineno}
\usepackage{algorithm2e}
\usepackage{natbib}

\newtheorem{thm}{Theorem}

\newtheorem{lem}{Lemma} 
\newdefinition{rmk}{Remark}
\newproof{pf}{Proof}
\newproof{pot}{Proof of Theorem \ref{thm2}}
\newtheorem{example}{Example}
\newdefinition{defin}{Definition}
\newdefinition{coro}{Corollary}

\begin{document}

\begin{frontmatter}

\title{A random free-boundary diffusive logistic model: Analysis, computing and simulation}

\author[IMM]{M.-C. Casab\'{a}n\corref{cor1}}
\ead{macabar@imm.upv.es}
\author[IMM]{R. Company}
\ead{rcompany@imm.upv.es}
\author[UC]{V.N. Egorova}
\ead{vera.egorova@unican.es}
\author[IMM]{L. J\'{o}dar}
\ead{ljodar@imm.upv.es}

\cortext[cor1]{Corresponding author. Tel.: +34 (96)3879144.}

\address[IMM]{Instituto Universitario de Matem\'{a}tica Multidisciplinar, Building 8G, access C, 2nd floor, Universitat Polit\`{e}cnica de Val\`{e}ncia, Camino de Vera s/n, 46022 Valencia, Spain}
\address[UC]{Depto. de Matemática Aplicada y Ciencias de la Computaci\'{o}n, Universidad de Cantabria, Avda. de los Castros, s/n, 39005 Santander, Spain}

\begin{abstract}
A free boundary diffusive logistic model finds application in many different fields from biological invasion to wildfire propagation. However, many of these processes show a random nature and contain uncertainties in the parameters. In this paper we extend the diffusive logistic model with unknown moving front to the random scenario by assuming that the involved parameters have a finite degree of randomness. The resulting mathematical model becomes a random free boundary partial differential problem and it is addressed numerically combining the finite difference method with two approaches for the treatment of the moving front. Firstly, we propose a front-fixing transformation, reshaping the original random free boundary domain into a fixed deterministic one. A second approach is using the front-tracking method to capture the evolution of the moving front adapted to the random framework. Statistical moments of the approximating solution stochastic process and the stochastic moving boundary solution are calculated by the Monte Carlo technique. Qualitative numerical analysis establishes the stability and positivity conditions. Numerical examples are provided to compare both approaches, study the spreading-vanishing dichotomy, prove qualitative properties of the schemes and show the numerical convergence.
\end{abstract}


\end{frontmatter}


\section{Introduction and motivation}    \label{sec:Introduction}

Random partial differential equations (RPDEs) present a dynamic area of research which deals with the modelling of systems where some parameters are subject to random variations. The behaviour of the system is influenced by randomness, thus leading to solutions that are random processes themselves.  This model is relevant when studying diffusive logistic models of Fisher-KPP type \citep{Fisher_1937, Kolmogorov_1937}, which seek to capture the evolution of a population that both grows and spreads within a particular region. In these models, certain parameters, such as the growth rate or diffusion coefficient, are assumed to vary randomly becoming random variables (r.v.'s) and reflecting the inherent variability in natural phenomena.

Dynamic population models have found extensive application in diverse fields \cite{Brauer2012}. Ecologists use them to better understand and predict the spread of species, taking into account the stochastic nature of ecological systems \citep{Acevedo2012, Zhou2014}. Similarly, they are utilized in predicting the spread of wildfires, where the speed and direction of fire propagation are highly sensitive to random factors like wind and moisture content \citep{Pagnini_2014}. In biology, diffusive logistic models with random parameters help describe the dynamics of cell populations due to the cancer progression \citep{Durrett2010}.

Apart form the study of the distribution of the population in  space and time, a relevant aspect of the diffusive logistic models is the treatment of the spreading front of the population. The consideration of this moving front turns the problem into a free boundary one where the localization of the front is an additional unknown. Free boundary problems cover a wide framework of applications in science and engineering \cite{Muntean2023, Vynnycky2023,Miklavcic2023,Friedman2015,DAcunto2021,He2023,Lu2020,Chadam1993,
Vromans2018}. 
 The free boundary diffusive logistic model was introduced in \citep{Du2010Spreading-VanishingBoundary} by adding a Stefan's type condition for the speed of the moving front, see also the illustrative survey \cite{Du2022}. 
In the random framework the position of the free boundary depends on the interaction between the population's growth, its diffusion characteristics, and the randomly varying environmental parameters. This makes the free boundary a stochastic process (s.p.), similar to the population density, adding an extra layer of complexity to the models. In the present study we follow the stochastic approach based on the mean square (m.s.) calculus \citep{TTSoong} where  there is a known
uncertainty, that is, the s.p.'s and r.v.'s can follow a wide range of probability distributions such as beta, exponential, Gaussian, etc.

Given the complexity and the high-dimensional nature of these random problems, random numerical methods play a crucial role in their analysis \citep{CasabanCompanyJodar2021}. Such methods allow us to approximate the solutions of these models, providing insights into the behaviour of the system under various scenarios. They can handle the random parameters and the dynamic free boundary, offering a practical way to study these stochastic systems. As such, developing and refining numerical methods for RPDEs with free boundaries is of great importance in exploiting these models to their full potential. 

In this work we present random numerical methods for a random extension of the free boundary diffusive logistic model of Du and Lin  \citep{Du2011Spreading-vanishingII, Du2020}
  with radial symmetry.  The random  population density of a spreading species,  $u(\omega)=u(r,t;\, \omega)$,  and the random moving boundary, $H(t; \, \omega)$, are 2-stochastic processes (s.p.'s) defined in a complete probability space $(\Omega, \mathcal{F}, \mathbb{P})$. In accordance with a previous work about random moving boundary phase-change problems \citep{CasabanCompanyJodarMATCOM2023}, for the sake of practical application we assume that the uncertainty is limited to  $p$-degrees of randomness \citep[p.37]{TTSoong}, i.e. they depend on a finite number $p$ of random variables (r.v.'s):
\begin{eqnarray}  \label{eq:DegreeRandomness}
u(\omega) &    = & \tilde{F}(r,t; A_1(\omega),A_2(\omega),\ldots,A_p(\omega))\,,\\
H(t;\omega) & = & \tilde{G}(t;B_1(\omega),B_2(\omega),\cdots,B_p(\omega))\,,
\end{eqnarray}
where
\begin{eqnarray} 
A_i(\omega)\,, \ B_i(\omega)\,, \  i=1,\ldots,n,  \text{are mutually independent r.v.'s}; \\
\tilde{F} \  \text{is a second order differentiable real function on variable} \ r \ \text{and first order differentiable on variable} \ t;   \\
\tilde{G} \ \text{is a differentiable real function on variable} \ t.   \label{eq:DegreeRandomness_functions}
\end{eqnarray}

Under hypotheses \eqref{eq:DegreeRandomness}--\eqref{eq:DegreeRandomness_functions} the s.p.'s $u(\omega)$ and $H(t;\omega)$ have sample differentiable trajectories (realizations), i.e., for a fixed event, $\omega_{\ell}\in{\Omega}$, they are differentiable functions of the real variables $(r,t)$ and $t$, respectively. Additional hypotheses related to the mean square  (m.s.) calculus will be impose later.

Inspired in the seminal radially symmetry deterministic model \citep{Du2011Spreading-vanishingII}, this paper is focused on the following random  free boundary diffusive logistic model of Stefan type in  the m.s. sense:
 
\begin{eqnarray}
    u_t(\omega) = D(\omega)\left( u_{rr}(\omega) + \frac{1}{r}u_r(\omega)\right)+ u(\omega)\left( \alpha(r) - \beta(r)\, u(\omega)\right), \quad t>0, \; 0<r< H(t;\omega),\ \omega\in{\Omega}\,, \label{eq:PDE}\\
    H'(t;\omega) = -\eta(\omega) \ u_r(H(t;\omega),t;\omega), \quad t>0,\ \ \omega\in{\Omega}\,, \label{eq:stefan}
\end{eqnarray}
subject to the initial and boundary conditions
\begin{eqnarray}
H(0;\, \omega)&= &H_0,\;    u(r,0; \,  \omega) = u_0(r), \quad 0\leq r \leq H_0, \label{eq:IC}\\
   u_r (0,t;\, \omega) &=& 0, \quad u(H(t;\, \omega), t;\,\omega) = 0,\quad t>0. \label{eq:BC}
\end{eqnarray}
The positive random variable (r.v.) $D(\omega)>0$ in the RPDE  \eqref{eq:PDE} denotes the diffusion rate   
and it is bounded such that
\begin{equation}  \label{eq:bound}
0 < d_1 \leq D(\omega)  \leq d_2\,, \quad  \text{for every} \ \omega\in{\Omega}\,.
\end{equation}
The r.v. $\eta(\omega)$, in the Stefan condition \eqref{eq:stefan}, denotes a  positive r.v. meaning the proportionality between the random moving boundary speed, $H'(t;\omega)$, and the random population gradient at the front. It is assumed there exists a bound $\eta_0$ such that
\begin{equation} \label{eq:bound_eta}
0 < \eta_0 \leq \eta (\omega), \quad \text{for every} \, \omega\in{\Omega}.
\end{equation}

The function $\alpha(r)$ is the intrinsic growth rate and $\alpha(r)/\beta(r)$ is the habitat carrying capacity of the species. Both, $\alpha(r)$ and $\beta(r)$ are  bounded continuous real functions satisfying
\begin{equation}\label{eq:kappa}
    \exists \ \kappa_1, \kappa_2 >0: \quad \kappa_1 \leq \alpha(r), \, \beta(r) \leq \kappa_2, \quad \forall r \in[0,\infty).
\end{equation}
Moreover, initial population density  function $u_0(r)$ satisfies  
\begin{equation}  \label{eq:ConditionsCI}
    u_0 \in C^2\left([0, H_0]\right), \quad u_0'(0) = u_0(H_0) =0, \quad u_0(r)>0, \ \  \forall r\in[0,H_0)\,.
\end{equation}
 $H_0$ is the radius value for the circular region where the initial population is confined.


The aim of this paper is to construct stable numerical finite difference schemes (RFDS's) for the random diffusive logistic model \eqref{eq:PDE}--\eqref{eq:ConditionsCI} preserving the qualitative characteristics of its solution. Specifically, we focus on the study and the comparison of two methods: the random Front-Fixing (FF) method and the random Front-Tracking (FT) method.
The random FF method, which was previously investigated for the deterministic case in \citep{Casaban_deterministic}, is a transformation-based technique. In this method, we perform a change of variables that maps the moving boundary to a fixed one. This simplifies the problem, as we no longer need to track a moving boundary. However, the transformed problem often becomes more complex, and this method requires careful consideration of the transformation used. One of the defining traits of the random FF method is its ability to preserve the number of grid points, which contributes to its stability.
On the other hand, the random FT method takes a direct approach by tracking the moving boundary as part of the solution process. Instead of transforming the problem, this method adapts to the changing domain by dynamically updating the grid. This can provide a more accurate representation of the problem, especially when the boundary movement is substantial. Despite the computational challenges of dynamically handling the grid, the random FT method ensures consistent step-size in the original variables, which is essential for accuracy. Throughout the paper, we will discuss and compare these methods highlighting their strengths and weaknesses. We aim to provide a comprehensive understanding of their performance in solving the RPDE problem \eqref{eq:PDE}--\eqref{eq:ConditionsCI}.

The organization of the rest of this paper is as follows. Section \ref{sec:Preliminaries} recalls some definitions and concepts related to the mean square (m.s.) and mean fourth (m.f.) calculus, as well as some qualitative properties of the solution in the deterministic scenario. In Section \ref{sec:FFmethod}, we discuss the application of the random FF method in the context of RPDEs and the construction of a RFDS, detailing its formulation, algorithm, and numerical analysis. Section \ref{sec:FTmethod} focuses on the random FT method, including a detailed discussion on its mechanism, application nuances, and numerical characteristics. Section \ref{sec:MCmethod}, we explore the use of the Monte Carlo technique in order to avoid the computation drawbacks concerning with the storage information overload through the iterations of the RFDS. Subsequently, in Section \ref{sec:NumericalExamples}, we compare the random FF and FT methods, explore the spreading-vanishing dichotomy, and study numerical convergence through practical examples, which help evaluate the effectiveness of these methods. The paper concludes with a summary of our findings and their implications in Section \ref{sec:Conclusions}.

\section{Preliminaries}  \label{sec:Preliminaries}
This section is splitted into two parts. For the sake of clarity in the presentation, the first part is devoted to recall some definitions and concepts related to the mean square (m.s.) and mean four (m.f.) calculus, \citep{TTSoong}, \citep{VillafuerteBraumannCortesJodar2010}. In the second one we remember some qualitatives properties of the solution given in  \citep{Du2011Spreading-vanishingII}.
Furthermore we provide a result that  will be used for the random scenario.

\subsection{Random operational calculus}   \label{subsec:Preliminaries_msCalculus}
 A real r.v. $V(\omega)$ is called a $p$ random variable ($p$-r.v.), on the probability space $\left( \Omega ,\, \mathcal{F}, \, \mathcal{P}\right)$ if satisfies that
$ \mathbb{E} \left[ |V(\omega)|^p \right] < + \infty$,
where $\mathbb{E}[\cdot]$ denotes the expectation operator. We denote by $L_p(\Omega)$ the space of all the $p$-r.v.'s endowed with the norm
\[  \left\|  V(\omega) \right\|_{p} = \left( \mathbb{E}[|V(\omega)|^{p}] \right)^{1/p} = \left( \int_{\Omega} \left| V(\omega) \right|^p \, f _V(\omega) \, \mathrm{d}\omega  \right)^{1/p} < + \infty ,\]
where $f_V$ denotes the density function of the r.v. $V(\omega)$ and $\omega\in{\Omega}$ an event of the sample space $\Omega$.  The spaces $(L_p(\Omega),\| V(\omega) \|_p)$, $p>1$, are Banach spaces \citep[Chap. 4]{Okendal}, and they verify that $L_q(\Omega) \subset L_p(\Omega)$ for $q > p$. In particular, the cases  $p=2$ and $p=4$, that is, $L_2(\Omega)$ and $L_4(\Omega)$, represent the mean square (m.s.) and the mean four (m.f.) spaces, respectively. We will build the solution stochastic processes (s.p.'s) on these spaces.

The general treatment of the m.s. solutions of  RPDE requires some basic operational tools such as the derivative of the product of two s.p.'s.  or the chain rule for the composition of s.p. Difficulties with the $p$-mean operational calculus of the product arise due to the fact that the $p$-norm is not sub multiplicative, that is, it does not satisfy the Banach inequality $\|U(\omega)V(\omega) \|_p \leq \|U(\omega) \|_p \, \| V(\omega)\|_p$. 
In particular, m.s. differentiation product rule for m.s. differentiable processes holds when one of the factors is deterministic or independent. But these strong restrictions can be removed if both processes are m.f. differentiable, \citep{VillafuerteBraumannCortesJodar2010}. This fact links the m.s. calculus with the m.f. calculus. Then, a m.f. calculus is necessary to get the m.s. results for solving linear RPDE in the m.s. sense. By applying Schwarz inequality \citep[Page 43]{TTSoong}, one obtains the following inequality between $p$-norms:
\[ \|U(\omega)V(\omega) \|_p \leq \|U(\omega) \|_{2p} \, \| V(\omega)\|_{2p}\,,  \quad \forall U(\omega), V(\omega)\in{L_{2p}(\Omega)}, \]
that will be applied later for $p = 2$.

Let be $T \subset \mathbb{R}$, a subsequence of $t$-indexed r.v.'s $\{  V(t;\omega) \, : \, t\in{T} \}$, is called a stochastic process (s.p.) of order $p$ ($p$-s.p.) if for each $t_f \in{T}$ fixed, the r.v. $V (t_f ;\omega) \in{L_p(\Omega)}$.
In $L_p(\Omega)$ spaces, $p$-mean convergence is referred to as the corresponding $p$-norm, when $p = 2$, it is called m.s. convergence, for $p = 4$, m.f. convergence and so on. We recall that $q$-mean convergence entails $p$-convergence, whenever $q > p$. A $p$-s.p. $V(t;\omega)$ is called $p$-mean continuous at $t\in{T}$ if
\[ \left\| V(t+\delta;\omega)-V(t;\omega) \right\|_p \rightarrow 0 \quad \text{as} \ \delta \rightarrow 0, \ \ t,t+\delta\in{T}\,.\]
A $p$-s.p. is called $p$-mean differentiable at $t\in{T}$ if there exists the $p$-derivative of $V(t;\omega)$, namely  $V^{\prime}(t;\omega)$, verifying
\[ \left\| \frac{V(t+\delta;\omega)-V(t;\omega)}{\delta} - V^{\prime}(t;\omega) \right\|_p \rightarrow 0 \quad \text{as} \ \delta \rightarrow 0, \ \ t,t+\delta\in{T}\,.\]

\subsection{Spreading-Vanishing dichotomy of the model}   \label{subsec:Dichotomy-Preliminars}
In the deterministic case the authors in   \citep{Du2010Spreading-VanishingBoundary, Du2011Spreading-vanishingII} showed that a spreading-vanishing dichotomy  occurs for the unknown population depending on the initial value of the front $H_0$ and the parameter $\eta$. This result is suitable to any realization $\omega_{\ell}\in{\Omega}$ under hypotheses \eqref{eq:DegreeRandomness}--\eqref{eq:DegreeRandomness_functions}. In this way, Theorems 2.4, 2.5 and 2.10 of \citep{Du2011Spreading-vanishingII} allow to present the following criteria:

\begin{thm}[Spreading-Vanishing dichotomy \citep{Du2011Spreading-vanishingII}] \label{theo:SV}
	Let $\left\{ u(r,t;\omega_{\ell}), \, H(t;\omega_{\ell}) \right\}$ be the solution of the sample problem \eqref{eq:PDE}--\eqref{eq:BC} for a given realization $\omega_{\ell}\in{\Omega}$, then there is a positive constant $R^*(\omega_{\ell})$, such that 
	\begin{enumerate}
		\item If $H_0 \geq R^*(\omega_{\ell})$, then the population spreads;
		\item If $H_0< R^*(\omega_{\ell})$, then there exists ${\eta^*(\omega_{\ell})} >0$  depending on $u_0$ such that
		\begin{itemize}
			\item If $\eta(\omega_{\ell}) > \eta^*({\omega_{\ell}})$, then the population spreads,
			\item If $\eta(\omega_{\ell})  \leq \eta^*({\omega_{\ell}})$, then the population vanishes.
		\end{itemize}
	\end{enumerate}
\end{thm}
The threshold $R^*(\omega_{\ell})>0$ is the unique real value, see \citep[Theorem 2.1]{Du2011Spreading-vanishingII}, such that 
\begin{equation}   
	\lambda_1 (D(\omega_{\ell}), \alpha(r), R^*(\omega_\ell)) =1, 
\end{equation}
where $\lambda_1 (D(\omega_{\ell}), \alpha(r), R^*(\omega_{\ell}))$ is the principal eigenvalue of the  problem
\begin{equation}
	\begin{cases}
		D(\omega_{\ell}) \, \left( \frac{d^2\, \phi}{dr^2} + \frac{1}{r}\frac{d\, \phi}{dr} \right)
		+ \lambda(\omega_\ell) \, \alpha(r) \, \phi = 0, \quad 0<r<R^*(\omega_{\ell}),\\
		\phi\left(R^*(\omega_{\ell})\right) =0. 
	\end{cases}
	\label{eq:PVF}
\end{equation}
In paper \citep{Casaban_deterministic}, the boundary value problem \eqref{eq:PVF} for $\lambda(\omega_\ell)=1$ is proved to be equivalent to the initial value problem
\begin{equation}
		\begin{cases}
			D(\omega_{\ell}) \, \left( \frac{d^2\, \phi}{dr^2} + \frac{1}{r}\frac{d\, \phi}{dr} \right) \phi + \alpha(r) \phi = 0, \quad r>0,\\
			\phi(0) = C,\\
			\phi'(0) = 0,
		\end{cases}
		\label{eq:IVP}
	\end{equation}
where $C>0$ is an arbitrary constant. Then, $R^*(\omega_\ell)$ is the first positive root of  $\phi(r_\ell;\omega_\ell) =0$. 
We are interested in finding an upper bound of $\left\{R^*(\omega_\ell) \, : \, \omega_{\ell}\in{\Omega}\right\}$, that guarantees spreading behaviour for all realizations when the initial front, $H_0$, exceeds this bound, that is, case 1 of Theorem \ref{theo:SV}. For this propose we establish the following result on inequalities for a particular second order linear ordinary differential initial value problem, see \citep{Azbelev}.

\begin{lem}  \label{lemma:R*}
Let $y_i(r)$, $i=1,2$, be the solutions of the following ordinary differential equations (ODE's) problems:
\begin{eqnarray}  \label{eq:LemmaProblemR*}
\left.\begin{array}{ll}
y_i^{\prime \prime}(r)+ \frac{1}{r}\, y_i^{\prime}(r) + q_i(r) \, y_i(r) = 0 & r>0\,, \\
  y_i(0)=C, \quad y_i^{\prime}(0)=0\,,   &   C>0 \ \text{(constant)}
\end{array} \right\}  \,, \  i=1,2,
\end{eqnarray}
being $q_i(r)$ a positive continuous function. Let us consider that 
\begin{equation} \label{eq:Lemma_q}
q_1(r) \leq q_2(r) \,, \quad  r>0\,,
\end{equation}
 and let $R_i$, $i=1,2$, be the first root of the equation $y_1(r)=0$, and $y_2(r)=0$, respectively. Then, it holds that
\begin{equation} \label{eq:LemmaRBound}
R_1  \geq R_2\,.
\end{equation}
\end{lem}

\begin{pf}
Let us define $y(r) = y_1(r) - y_2(r)$, then subtracting the equations given in \eqref{eq:LemmaProblemR*}  one gets
\begin{equation}   \label{eq:Lemma_eq1}
y^{\prime \prime}(r)+ \frac{1}{r}\, y^{\prime}(r) + q_1(r)\,  y(r) =
\left( q_2(r) - q_1(r)   \right) \, y_2(r)\,.
\end{equation}
Let us take $R=\min \left( R_1,\ R_2 \right)$. By the continuity of the solution $y_2(r)$  we have
\begin{equation} \label{eq:Lemma_eq2}
y_2(r) >0\,,  \quad 0 \leq r < R\,. 
\end{equation}
  From \eqref{eq:LemmaProblemR*}--\eqref{eq:Lemma_q} , \eqref{eq:Lemma_eq1}--\eqref{eq:Lemma_eq2} the solution 
 $y(r)=y_1(r)-y_2(r)$ satisfies the inequality problem
\begin{eqnarray}  \label{eq:Lemma_eq3}
\left.\begin{array}{ll}
y^{\prime \prime}(r)+ \frac{1}{r}\, y^{\prime}(r) + q_1(r) \, y(r) \geq 0\,, & 0 < r < R\,, \\
  y(0)= y^{\prime}(0)=0\,,   &   
\end{array} \right\}  
\end{eqnarray}
Further, the unique solution of the following problem
\begin{eqnarray}  \label{eq:Lemma_eq4}
\left.\begin{array}{ll}
z^{\prime \prime}(r)+ \frac{1}{r}\, z^{\prime}(r) + q_1(r) \, z(r) = 0\,, & 0 < r < R\,, \\
 z(0)= z^{\prime}(0)=0\,,   &   
\end{array} \right\}  
\end{eqnarray}
is the null function, $z(r)=0(r)$. It can be proved that satisfies $y(r) \geq z(r)=0(r)$ in  $0 <r < R$. 
 In fact, subtracting ODE in \eqref{eq:Lemma_eq4} from inequality in \eqref{eq:Lemma_eq3}, and denoting 
 \[ f(r)=\frac{1}{r}+ \frac{q_1(r)\,\left( y(r) - z(r)  \right)}{y^{\prime}(r)-z^{\prime}(r) } \,, \]
 one gets
 \begin{equation}  \label{eq:Lemma_eq5}
y^{\prime \prime}(r) - z^{\prime \prime}(r)+ f(r) \, \left(
y^{\prime}(r) - z^{\prime}(r) \right)  \geq  0\,, \quad  0 < r < R\,.
\end{equation}
By multiplying both sides of \eqref{eq:Lemma_eq5} by  $\mathrm{exp}\left(\int_0^{\, r} \, f(s) \, \mathrm{d}s\right)$ one gets
\begin{equation}  \label{eq:Lemma_eq6}
\frac{d}{dr}   \left[ \left( y^{\prime }(r) - z^{\prime }(r) \right)  \, \mathrm{exp} \left( \int_0^{\, r} \, f(s;\omega_\ell) \, \mathrm{d}s \right) \right] \ \  \geq 0\,.
\end{equation}
Integrating \eqref{eq:Lemma_eq6} between $0$ and $\bar{r}$ with $0 < \bar{r} < R$, we have
\begin{equation}  \label{eq:Lemma_eq7}
\left( y^{\prime}(\bar{r}) - z^{\prime}(\bar{r}) \right) \,
 \mathrm{exp} \left( \int_0^{\,\bar{r}} \, f(s) \, \mathrm{d}s \right)  - 
 \left( y^{\prime}(0) - z^{\prime}(0) \right) \   \geq  0\,.
\end{equation}
Renaming $r = \bar{r}$ and taking into account the positivity of the exponential in \eqref{eq:Lemma_eq7} one gets
\begin{equation}   \label{eq:Lemma_eq8}
 y^{\prime}(r) \geq z^{\prime}(r )=0(r)\,, \quad 0 <r < R\,.
\end{equation}
Integrating in \eqref{eq:Lemma_eq8}  and taking into account initial condition of \eqref{eq:Lemma_eq3} we have
\begin{equation}
y(r)  \geq 0\,, \quad 0 < r< R\,.
\end{equation}
Hence by definition of $y(r)$ one gets
\[  y_1(r) \geq  y_2(r) \,, \quad 0 < r < R\,, \]
consequently  the result  \eqref{eq:LemmaRBound} is proved.  $\qquad \square$
\end{pf}

Let us recall the boundedness condition \eqref{eq:bound} for the r.v. $D(\omega)$ and consider $R^*(d_2)$,  the first positive root of  the solution $\phi(r;d_2)$ of ODE problem  
\begin{equation} 
		\begin{cases}
			d_2 \, \left( \frac{d^2\, \phi}{dr^2} + \frac{1}{r}\frac{d\, \phi}{dr} \right) \phi + \alpha(r) \phi = 0, \quad r>0, \\
			\phi(0) = C,\\
			\phi'(0) = 0,
		\end{cases}
		\label{eq:IVP}
	\end{equation}
Let us take  $q_1(r)$ and $q_2(r)$ defined by
\[q_1(r)=\frac{\alpha(r)}{d_2}\,, \quad q_2(r)=q(r;\omega_\ell) =\frac{\alpha(r)}{D(\omega_\ell)}, \ \ \omega_\ell\in{\Omega}\,.\] 
Then from \eqref{eq:bound} one gets $q_1(r) \leq q_2(r)$ and from Lemma \ref{lemma:R*} 
\[
R^*(d_2) \geq R^*(\omega_\ell)\,, \quad \text{for a fixed}\ \omega_\ell\in{\Omega}\,. \]
 Consequently, from Theorem \ref{theo:SV} the spreading of the population for all realizations $\omega\in{\Omega}$ is guaranteed  if the initial front, $H_0$, satisfies
\begin{equation}  \label{eq:ConditionBoundH0}
H_0 \geq R^*(d_2) \,.
\end{equation}

\section{A random Front-Fixing (FF) method}  \label{sec:FFmethod}
In this section we present a random FF method together with the construction of a RFDS for solve numerically the radial symmetric diffusive logistic  model \eqref{eq:PDE}--\eqref{eq:BC}. Dealing with populations the positivity of the numerical solution is an important issue and it is treated together with the stability analysis in the last part of the section.

The FF method is  an immobilization technique for the free boundary by using Landau-type transformations. For the random problem  \eqref{eq:PDE}--\eqref{eq:BC} we use the change of variables 
\begin{equation}\label{eq:landau}
    z = \frac{r}{H(t;\, \omega)}, \quad v(z,t;\, \omega) = u(r,t;\, \omega)\,.
\end{equation}
Note that the new spatial variable $z$ is a fixed proportion  $0 \leq z \leq 1$ for each $r\in{[0,\,H(t;\, \omega) ]}$ at fixed time $t>0$.  From this point of view, the original spatial variable $r=z \, H(t;\omega)$ depends on every event $\omega\in{\Omega}$ and  it becomes a r.v. for every fixed $t$. Moreover, deterministic parameters $\alpha(r)$ and $\beta(r)$ become random ones as well
\begin{eqnarray}
    \alpha(r) =  \alpha(z\, H(t;\, \omega)) =a(z,t;\, \omega) ,\\
   \beta(r)  = \beta(z\, H(t;\, \omega)) =  b(z,t;\, \omega) ,
\end{eqnarray}
with the same constraint 
\begin{equation}
    0< \kappa_1 \leq a(z,t;\, \omega), \; b(z,t;\, \omega) \leq \kappa_2\,, \qquad t>0, \; 0\leq z \leq 1, \; \omega\in \Omega.
\end{equation}

By denoting $G(t;\, \omega) = H^2(t;\, \omega)$ and substituting \eqref{eq:landau} into the problem \eqref{eq:PDE}--\eqref{eq:BC}, one gets

\begin{multline}
    G(t;\, \omega) \,v_t  = D(\omega) \,  v_{zz} + \left( \frac{D(\omega)}{z} + \frac{z}{2}G'(t;\, \omega)\right) \, v_z 
  \\  + G(t;\, \omega) \left[ a(z,t;\, \omega) - b(z,t;\, \omega) \, v \right] \ v, 
 \qquad    \qquad  t>0, \; 0<z< 1, \omega\in \Omega,\label{eq:PDE_tr}
\end{multline}
\begin{equation}
G'(t;\, \omega)  =  -2\, \eta(\omega) \, v_z(1,t;\, \omega), \quad t>0, \; \omega\in\Omega,\label{eq:stefan_tr}
\end{equation}
subject to the initial and boundary conditions
\begin{eqnarray}
G(0)= H_0^2, &  v(z,0) = u_0(z\, H_0^2),  & 0\leq z \leq 1,  \label{eq:IC_tr}\\
   v_z (0,t;\, \omega) = 0, & v(1, t;\, \omega) = 0, &  t>0,  \, \omega \in \Omega.\label{eq:BC_tr}
\end{eqnarray}

\begin{rmk}
The m.s. operational calculus developed in \eqref{eq:PDE_tr} and \eqref{eq:stefan_tr} is legitimised by imposing that the s.p.'s:  $v^2 $, $v_z$, $v_t$, $G(t;\omega)$, $G^{\prime}(t;\omega)$ lie in $L_4(\Omega)$, see Preliminaries in Section \ref{subsec:Preliminaries_msCalculus}.
\end{rmk}

\subsection{Construction of the random finite difference scheme for FF method (RFDS-FF)}

We consider the transformed problem \eqref{eq:PDE_tr}--\eqref{eq:BC_tr}. 
Let us define a uniform grid $z_j = jh$, $h = 1/M$, $t^n = nk$, $k = T/N$, where $N$ and $M$ are two given positive integer numbers and $T$ is a fixed time horizon. We denote by $v_j^n(\omega) \approx v(z_j,t^n;\omega)$ the numerical approximation of the solution s.p. $v(z,t;\omega)$, $\omega\in{\Omega}$, in a mesh point $(z_j,t^n)$. In addition we introduce  the following notation
\begin{equation}
g^n(\omega) \approx G(t^n;\omega), \quad a_j^n(\omega) = a(z_j, t^n;\omega), \quad b_j^n(\omega) = b(z_j,t^n;\omega)\,. 
\end{equation}
By using a forward first-order approximation for the time derivatives and a centred second-order approximations for the spatial derivatives in 
\eqref{eq:PDE_tr} one gets the following  explicit RFDS-FF:
\begin{equation}\label{eq:scheme_coef}
    v_j^{n+1}(\omega) = A_j^n(\omega)\, v_{j-1}^n(\omega) + B_j^n(\omega) \,v_j^n(\omega) +
     C_j^n(\omega) \, v_{j+1}^n(\omega)\,, \quad 1 \leq j \leq M-1\,, \ 0 \leq n \leq N-1\,, \ \omega\in{\Omega}\,,
\end{equation}
where
\begin{eqnarray} \label{eq:Ajn-Bjn-Cjn}
\left. \begin{array}{rcl}
    A_j^n(\omega) &= &\dfrac{D (\omega)\, k}{h^2\,g^n(\omega)} - \dfrac{D(\omega) k}{2h \,g^n(\omega) \,z_j} - \dfrac{z_j}{4h}\left(\dfrac{g^{n+1}(\omega)}{g^n(\omega)}-1\right) \\
    \\
    B_j^n(\omega) & = & 1+ k \left(a_j^n(\omega) - b_j^n (\omega)\, v_j^n(\omega) - \dfrac{2D(\omega)\, k}{h^2 g^n(\omega)}\right)  \\
    \\
       C_j^n(\omega) &= &\dfrac{D(\omega) \, k}{h^2\,g^n(\omega)} + \dfrac{D(\omega) \, k}{2h\, g^n(\omega)\, z_j} + \dfrac{z_j}{4h}\left(\dfrac{g^{n+1}(\omega)}{g^n(\omega)}-1\right) 
 \end{array} \right\} \quad 0 \leq n \leq N-1\,, \ \omega\in{\Omega} \,.
\end{eqnarray}
The time derivative in the Stefan condition \eqref{eq:stefan_tr} is discretized by using  a forward first-order approximation while the spatial derivative at the fixed front  $z=1$ is discretized by the second-order left side finite difference approximation obtaining
\begin{equation} \label{eq:scheme_g}
    \frac{g^{n+1}(\omega)-g^n(\omega)}{k} = -\frac{\eta(\omega)}{h} \, \left( 3v_M^n(\omega) -4v_{M-1}^n(\omega)+ v_{M-2}^n(\omega) \right)\,, \quad 0 \leq n \leq N-1\,.
\end{equation}
From \eqref{eq:scheme_g} and using \eqref{eq:BC_tr} the explicit random scheme for the evolution of the random moving front is given by
\begin{equation}  \label{eq:scheme_g2}
 g^{n+1}(\omega) = g^n(\omega) + \dfrac{k}{h}\eta(\omega)\, (4v_{M-1}^n(\omega)-v_{M-2}^n(\omega))\,.
 \end{equation}
Boundary conditions \eqref{eq:BC_tr} are discretized as
\begin{equation}\label{eq:scheme_BC}
    v_z(0,t^n;\omega) \approx \frac{-3v_0^n(\omega) + 4v_1^n(\omega) -v_2^n(\omega)}{2h} = 0, \quad v_M^n(\omega) = 0\,, \quad 0 \leq n \leq N\,.
\end{equation}
Once numerical solution for the transformed RPDE is computed, the numerical solution of the original Stefan problem at $t=T$ is  found by the inverse transformation
\begin{equation}
    r_j \approx z_j \, H(T;\, \omega), \quad u(r_j,T; \, \omega) \approx v_j^{N}(\omega). 
\end{equation}

\subsection{Stability and qualitative properties of the numerical solutions s.p.'s}   \label{subsec:FFmethod_properties}
In this section we show qualitative properties of the RFDS-FF \eqref{eq:scheme_coef}--\eqref{eq:Ajn-Bjn-Cjn}, \eqref{eq:scheme_g2}--\eqref{eq:scheme_BC}. In particular, we will proof that the RFDS-FF is conditionally stable and preserves the positivity and the monotonicity of the random numerical s.p.'s.
For the sake of clarity in the presentation, we start recalling some definitions \citep{CasabanCompanyJodarMATCOM2023}.

\begin{defin}    \label{def:non_increasing_monotone}
    The  numerical solution s.p. $\left\lbrace v_j^n(\omega)\right \rbrace$, $\omega \in \Omega$, of a random RFDS is said to be non-increasing monotone in the spatial index $j = 0,\ldots, M-1$, if
    \begin{equation} 
        v_j^{n}(\omega_{\ell}) \geq v_{j+1}^n(\omega_{\ell}), \quad  0 \leq j \leq M-1, \; 0\leq n \leq  N,
    \end{equation}
for every $\omega_{\ell}\in \Omega$.
\end{defin}

\begin{defin}    \label{def:non_increasing}
    The  numerical free boundary s.p. $\left\lbrace g^n(\omega)\right \rbrace$, $\omega \in \Omega$ of a random RFDS is said to be strictly increasing, if
    \begin{equation}
        g^{n}(\omega_{\ell}) < g^{n+1}(\omega_{\ell}), \quad  0\leq n \leq  N-1,
    \end{equation}
for every $\omega_{\ell}\in \Omega$.
\end{defin}

\begin{defin} \label{def:stability}
   A random numerical scheme is said to be ${\| \cdot \|}_{p}$-stable in the fixed station sense in $[0,\, 1]\times [0, \, T]$, if for every partition with $k=T/N$ and $h = 1/M$,  it is fulfilled that 
    \begin{equation}
      \left\|v_j^n(\omega_{\ell}) \right\|_p \leq C, \qquad  0\leq j\leq M, \; 0\leq n \leq N,
    \end{equation}
    where the constant $C$ is independent of the step-sizes $h$ and $k$  the time level $n$, for every realization $\omega_{\ell} \in \Omega$.
\end{defin}

In \citep{Casaban_deterministic} the free boundary diffusive logistic model is studied in the deterministic scenario proposing the Front-Fixing method combined with a finite difference scheme. Numerical analysis  developed provides qualitative properties of the numerical solution summarized in Theorems 2-6 of \citep{Casaban_deterministic}. In the present random scenario, we take into account that $D(\omega)>0$ is a bounded r.v. verifying condition $D(\omega)<d_2$ for every $\omega\in{\Omega}$, see \eqref{eq:bound}. As Definitions \ref{def:non_increasing_monotone}--\ref{def:stability} above are related to  the events of the numerical solution s.p.'s, Theorems 2-6 in \citep{Casaban_deterministic} provide the following result.
\begin{thm}  \label{theorem:properties}
    With the previous notation for small enough spatial step-size $h$, the random numerical solution $\left\lbrace v_j^n(\omega)\right \rbrace$,  and the random free-boundary $\left\lbrace g^n(\omega)\right\rbrace$,  $\omega\in \Omega$, satisfy
    \begin{enumerate}
        \item $v_j^n(\omega)>0$ for $0\leq j \leq M$, $0 \leq n\leq N$;
        \item $\left\lbrace v_j^n(\omega)\right \rbrace$ is  non-increasing monotone in the spatial index $j$;
        \item  $\left\lbrace g^n(\omega)\right\rbrace$ is strictly increasing;
        \item $\left\lbrace v_j^n(\omega)\right \rbrace$ is ${\| \cdot \|}_p$-stable in the fixed station sense;
    \end{enumerate}
    under the time step-size condition: 
    \begin{equation} \label{eq:ConditionStability_k}
        k < Q h^2,
    \end{equation}
    with
    \begin{equation}  \label{eq:Qi}
        Q = \min_{1\leq i\leq 3} Q_i,
    \end{equation}
\begin{equation}
    Q_1  =   \frac{H_0^2}{2d_2 + h^2 \alpha_1 H_0^2\left( \frac{\alpha_2\,\beta_2}{\alpha_1\,\beta_1} -1\right)  };
    \quad
    Q_2  = 
     \frac{H_0^2}{2d_2+ h^2\, \beta_2 \, H_0^2  (2M_0 - C_m)};
    \quad
    Q_3  = 
    \frac{4H_0^2}{9d_2 + 8h^2 \, \beta_2\, H_0^2 \,P_0 };
\end{equation}
where
\begin{equation}   \label{eq:Def_elem_Theorem2}
 C_0 = \sup_{r \in \mathbb{R}^+} \left\lbrace  \frac{\alpha(r)}{\beta(r)}\right\rbrace, \quad
    C_m = \inf_{r \in \mathbb{R}^+} \left\lbrace  \frac{\alpha(r)}{\beta(r)}\right\rbrace, \quad M_0 = \max_{0\leq r \leq H_0}\{u_0(r)\}, \quad P_0 = \max \left\lbrace M_0, \, C_0 \right\rbrace,
\end{equation}
and $\alpha_i$, $\beta_i$, $i=1,2$, verify
\begin{equation}\label{eq:kappa_bound}
   0 <  \kappa_1 \leq \alpha_1 \leq \alpha(r) \leq  \alpha_2 \leq \kappa_2, \qquad 
    0 < \kappa_1 \leq \beta_1 \leq \beta(r) \leq \beta_2 \leq \kappa_2,   \quad \forall r \in[0,\infty)\,.
\end{equation}
 \end{thm}

\subsection{Algorithm for computing the numerical  solution s.p.'s for FF method}   \label{subsec:Algorithm_FF}
Algorithm \ref{algo:EFDM} summarize the procedure of the numerical solution s.p.'s for every realization $\omega_\ell \in \Omega$.
 \RestyleAlgo{ruled}
\begin{algorithm}
\caption{Explicit RFDS-FF for the Stefan problem \eqref{eq:PDE}--\eqref{eq:BC} for a fixed $\omega_\ell\in \Omega$}\label{algo:EFDM}
\KwData{$D(\omega_\ell)$, $\eta(\omega_\ell)$, $\alpha(r)$, $\beta(r)$, $T$, $H_0$, $u_0$, $M$, $N$}
\KwResult{$r$, $u(r,T;\omega_\ell)$, $H(T;\omega_\ell)$}
$h \gets 1/M$, 

$k \gets T/N$ \text{verifying the stability condition}  \eqref{eq:ConditionStability_k}

$g_0 \gets H_0^2$

\For{$j = 0, \ldots, M$}{
$z_j \gets jh$, \
$r_j \gets H_0\, z_j$, \
$v^0_j = u_0(r_j)$
}
\For{$n = 0, \ldots, N-1$}{
    $g^{n+1}(\omega_\ell) \gets g^n(\omega_\ell) + \dfrac{k}{h}\eta(\omega_\ell)\, (4v_{M-1}^n(\omega_\ell)-v_{M-2}^n(\omega_\ell))$ \ \text{by}  \eqref{eq:scheme_g2}

     $H^{n+1}(\omega_\ell) \gets \sqrt{g^{n+1}(\omega_\ell)}$

    \For{ $j = 1,\ldots, M-1$}{
\text{Compute coefficients} $A_j^n(\omega_\ell)$, $B_j^n(\omega_\ell)$ and $C_j^n(\omega_\ell)$ \text{by} \eqref{eq:Ajn-Bjn-Cjn}

$v_j^{n+1}(\omega_\ell) \gets  A_j^n(\omega_\ell) \, v_{j-1}^n(\omega_\ell) + B_j^n(\omega_\ell) \, v_j^n(\omega_\ell) + C_j^n(\omega_\ell) \, v_{j+1}^n(\omega_\ell)$ \ \text{by}  \eqref{eq:scheme_coef}.
    }
    $v_0^{n+1}(\omega_\ell) \gets \left(4v_1^{n+1}(\omega_\ell) - v_2^{n+1}(\omega_\ell)\right) / 3$, \quad 
     $v_M^{n+1} (\omega_\ell)\gets 0$ \ \text{by} \eqref{eq:scheme_BC}.
  }
\For{$j = 0, \ldots, M$}{
$r_j \gets H^{N}(\omega_\ell) \, z_j$

$u_j (\omega_\ell)= v_j^{N}(\omega_\ell)$
}
\end{algorithm}


\section{A random Front-Tracking (FT) method }    \label{sec:FTmethod}

FF method and  Front-Tracking (FT) methods, \citep[Chap. 4]{Crank}, are usual techniques for free boundary problems of Stefan type in the deterministic framework. In this section we propose a random FT method for the studied random free boundary problem  \eqref{eq:PDE}--\eqref{eq:BC}. Comparison between random FF and random FT methods will be illustrated in Section \ref{subsection:comparisonFF-FT}. In FT methods the grid is computed for each time step due to the movement of the free boundary. They were classified by \citep{Gupta} as fixed and variable grid methods. We use the fixed grid method where the space-time domain is subdivided into a finite number of uniformly distributed nodes and the position of the boundary does not necessarily coincide with the mesh points.

\begin{figure}[h]
\begin{center}
	\includegraphics[width=10cm]{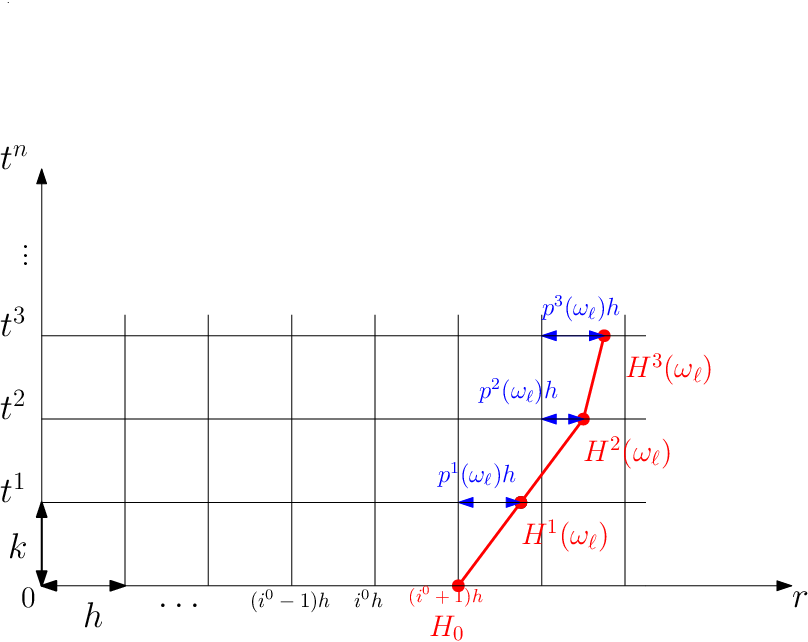}
\caption{ Evolution of the position of the moving boundary with FT method for a fixed realization $\omega_\ell\in{\Omega}$. The moving front, $H^n(\omega_\ell)$, in a time level $n$ is highlighted in red 
  and the fractional step forward in the space, $p^n(\omega_\ell)$, in a time level $n$ is marked in blue.}
\label{fig:FTracking_mesh}
\end{center}
\end{figure}

We consider the numerical domain $[0,L] \times [0,T]$, with the grid points $r_j=j\,h$, $t^n=n\,k$, where $h$ and $k$ are the space and time increments, respectively.  Figure \ref{fig:FTracking_mesh} shows the position of the moving boundary for a fixed realization $\omega_\ell\in{\Omega}$. The step sizes $(h,k)$ are fixed and it appears a fractional distance between the last interior mesh point and the localization of the moving boundary.
The numerical approximation of the solution s.p. $u(r_j,t^n;\omega)$ is denoted by $  u_j^n(\omega)$ and 
the approximation of the moving front s.p. $H(t^n;\omega)$ is denoted by  $H^n(\omega)$.
Let us denote $i^n(\omega)$ the spatial index of the last interior point of the non-zero population region,  that is, 
\[r_{i^n(\omega)} < H^n(\omega) \leq r_{i^n(\omega)+1}, \quad \omega\in{\Omega}. \]
 Then, for any realization $\omega\in{\Omega}$, there exists a distance parameter in each time level $n$ namely $p^n(\omega)$, such that
 \begin{equation}   \label{eq:_DefHn}
 H^n(\omega)=(i^n(\omega)+p^n(\omega))\,h\,, \quad 0 <p^n(\omega)\leq 1\,,  \ \ 0 \leq n \leq N.
\end{equation} 
From \eqref{eq:IC} one gets $H^0(\omega)=H_0$ and for initialization we take $p^0(\omega)=p^0=1$, $i^0(\omega)=i^0=M-1$ and $H^0(\omega)=H_0=(i^0+1)\,h=Mh$, for all realizations $\omega\in{\Omega}$.
Figure \ref{fig:FTracking_mesh} shows a simulation of the movement  of the front in the first time steps for a fixed realization $\omega_\ell\in{\Omega}$.

For the interior spatial points up to the last but one $r_j=jh$,  $1 \leq j \leq i^{n}(\omega)-1$,
we consider the following m.s. approximations
\begin{eqnarray}
u_t(r_j,t^n;\omega) & \approx & \frac{u_j^{n+1}(\omega)-u_j^{n}(\omega)}{k} \,, \quad \ 0 \leq n \leq N, \label{eq:u_tFT}\\
u_{rr}(r_j,t^n;\omega) & \approx & \frac{u_{j-1}^{n}(\omega)-2u_j^{n}(\omega)+u_{j+1}^n(\omega)}{h^2} \,, \quad \ 0 \leq n \leq N, \label{eq:u_rrFT}\\
u_r(r_j,t^n;\omega) & \approx & \frac{u_{j+1}^{n}(\omega)-u_{j-1}^{n}(\omega)}{2h}\,, \quad \ 0 \leq n \leq N \,.  \label{eq:u_rFT}
\end{eqnarray}
For $j=0$, taking into account the boundary \eqref{eq:BC}, $u_r(0,t^n;\omega)=0$, and using the second-order right-hand approximation for the first m.s. derivative in $r$ we have
\begin{equation}\label{eq:scheme_FT_j0}
   u_0^n(\omega)=\frac{4u_1^n(\omega)-u_2^n(\omega)}{3}\,, \ 1 \leq n \leq N.
\end{equation}
In order to approximate the m.s. spatial partial derivatives of the EDP \eqref{eq:PDE}  at the last interior point $r_{i^n(\omega)}=i^n(\omega)\, h$, a Lagrange interpolation passing through the points $r_{i^{n}(\omega)-1}$, $r_{i^n(\omega)}$ and $H^n(\omega)$ is used, see \citep[Pages 164-165]{Crank}, obtaining the approximations 
 \begin{eqnarray}
 \frac{\partial^2 u}{\partial r^2} \left( r_{i^n(\omega)},t^n;\omega \right)  & \approx & \frac{2}{h^2} \left( \frac{1}{p^n(\omega)+1} u_{i^{n}(\omega)-1}^n(\omega) -\frac{1}{p^n(\omega)} u_{i^{n}(\omega)}^n(\omega)  \right) \,,   \label{eq:Der_u_r2_FT}  \\
  \frac{\partial u}{\partial r} \left( r_{i^n(\omega)},t^n;\omega \right)  & \approx & \frac{1}{h} \left( - \frac{p^n(\omega)}{1+p^n(\omega)} u_{i^{n}(\omega)-1}^n(\omega) -\frac{1-p^n(\omega)}{p^n(\omega)} u_{i^{n}(\omega)}^n(\omega)     \right) \,,     \label{eq:Der_u_r_FT}
 \end{eqnarray}
 
In the boundary  $H^n(\omega)$, see \eqref{eq:_DefHn}, the same three point Lagrange interpolation formulae provides the following m.s. approximation for the first spatial m.s. derivative at the moving front involved in Stefan condition \eqref{eq:stefan}
\begin{eqnarray}
 \frac{\partial u}{\partial r} \left(H^{n}(\omega),t^n;\omega \right)  & \approx & \frac{1}{h} \left( \frac{p^n(\omega)}{p^n(\omega)+1} u_{i^{n}(\omega)-1}^n(\omega) -\frac{p^n(\omega)+1}{p^n(\omega)} u_{i^n(\omega)}^n(\omega)     \right) \,.   \label{eq:Der_u_r_H_FT} 
 \end{eqnarray}
The forward approximation of the time m.s. derivative in \eqref{eq:stefan}
\begin{equation}  \label{eq:DerStefan}
\frac{\mathrm{d}}{\mathrm{d}t}\left(H( t^n;\omega) \right) \approx \frac{H^{n+1}(\omega)-H^{n}(\omega)}{k} = \frac{\left( \Delta^n(\omega)-p^n(\omega) \right)\,h}{k}
\end{equation}
where  $\Delta^n(\omega) \, h$ denotes the distance between the last interior point in the time level $n$ and the localization of the moving front in time level $n+1$, that is
\begin{equation}   \label{eq:DefDelta_p}
\Delta^n(\omega) \, h= H^{n+1}(\omega) - i^{n}(\omega)\,h\,.
\end{equation}

From the approximations above the discretization of the RPDE \eqref{eq:PDE} takes the form
\begin{itemize}
\item[$\bullet$] For  $1 \leq j \leq i^{n}(\omega)-1$, using \eqref{eq:u_tFT}--\eqref{eq:u_rFT} one gets
\begin{equation} \label{eq:scheme_E1}
    u_j^{n+1}(\omega) = \tilde{A}_j^n(\omega) \, u_{j-1}^n(\omega) + \tilde{B}_j^n(\omega) \, u_j^n(\omega) + \tilde{C}_j^n(\omega) \, u_{j+1}^n(\omega),
\end{equation}
where
\begin{eqnarray} \label{eq:CoeficientsAjnBjn_Cjn_FT}
\left. \begin{array}{lcl}
    \tilde{A}_j^n(\omega) &= &\dfrac{k\, D(\omega)}{h^2} \left(1- \dfrac{1}{2j}  \right) \ > 0,  \\
    \\
    \tilde{B}_j^n(\omega) & = & 1+k\, \left[ - \dfrac{2D(\omega) }{h^2} + \alpha_j -\beta_j \, u_j^n(\omega)  \right],\quad \text{with}\ \alpha_j=\alpha(r_j), \ \beta_j=\beta(r_j)\,,\\
    \\
       \tilde{C}_j^n(\omega) &= &\dfrac{k\, D(\omega)}{h^2} \left(1+ \dfrac{1}{2j}  \right) > 0\,,
\end{array} \right\}
\end{eqnarray}

\item[$\bullet$] For the last interior point $r_{i^n(\omega)}=i^{n}(\omega)\,h$, we use the approximations \eqref{eq:u_tFT}, \eqref{eq:Der_u_r2_FT} and \eqref{eq:Der_u_r_FT} obtaining 
\begin{multline}  \label{eq:lastInteiorPoint}
u_{i^n(\omega)}^{n+1}(\omega)=u_{i^n(\omega)}^n(\omega) + k \, \left\{ \frac{D(\omega)}{h^2} \left[ \frac{u_{i^n(\omega)-1}^n(\omega)}{p^n(\omega)+1}\left( 2-\frac{p^n(\omega)}{i^n(\omega)}\right) - \frac{u_{i^n(\omega)}^n(\omega)}{p^n(\omega)}\left( 2+\frac{1-p^n(\omega)}{i^n(\omega)}\right)
 \right]  \right. \\
 \qquad \qquad  \left.+ u_{i^n(\omega)}^n(\omega)\left( \alpha_{i^n(\omega)}-\beta_{i^n(\omega)}\, u_{i^n(\omega)}^n(\omega)\right) \right\}\,.
\end{multline}  

\item[$\bullet$] For $j=0$
\begin{equation} \label{eq:u0_FT}
u_0^{n+1}(\omega) = \frac{4\, u_1^{n+1}(\omega)-u_2^{n+1}(\omega)}{3}
\end{equation}
\item[$\bullet$] For the advance of the front, given by the Stefan condition \eqref{eq:stefan},  from the approximations \eqref{eq:DerStefan}--\eqref{eq:DefDelta_p}
it  results that
\begin{equation}   \label{eq:p_FT}
\Delta^n(\omega) =p^n(\omega) + \frac{k \, \eta(\omega)}{h^2} \left(\frac{p^n(\omega)+1}{p^{n}(\omega)} \,u_{i^n(\omega)}^n(\omega)- \frac{p^n(\omega)}{p^{n}(\omega)+1}\,u_{i^n(\omega)-1}^n(\omega) \right) \,.
\end{equation}
Then the moving front at time level $n+1$ is given by 
\begin{equation}  \label{eq:Front_FT}
H^{n+1}(\omega)=\left( 	i^n(\omega) + \Delta^n(\omega) \right)\,h \,.
\end{equation}
\end{itemize}
Thus, the RFDS-FT for the problem \eqref{eq:PDE}--\eqref{eq:BC}     is given by  \eqref{eq:scheme_E1}--\eqref{eq:Front_FT}.

Dealing with the accuracy of the numerical scheme it is convenient to guarantee that the front advances up to one spatial-step $h$, that is,  $0 < H^{n+1}(\omega) -H^{n}(\omega) < h$. This means that the speed of advance of the random free front does not exceed the quotient $\frac{h}{k}$. Assuming that the speed of the moving front is expected to be non-increasing and from Stefan condition \eqref{eq:stefan} at $t=0$,  a suitable condition is the following
\[  \frac{h}{k} \geq  \, H^{\prime}(0)= - \eta(\omega) \, u^{\prime}_0(H_0)\,, \quad \omega\in{\Omega}. \]
From the bound \eqref{eq:bound_eta} one gets the time-step size constraint 
\begin{equation}   \label{eq:bound1_k_FT}
k \leq \frac{h}{\eta_0\, \left| u^{\prime}_0(H_0) \right|} \,.
\end{equation}
A relevant stability issue becomes from the fact that small values of the fractional spatial distance, $p^n(\omega)$, produces instabilities, see \eqref{eq:p_FT}. In accordance with \citep{Marshall_2} we  established a lower bound $0 < \epsilon < 1$ to avoid this undesirable situation. Then if in a time level $n$ one gets $0 < p^n(\omega) \leq \epsilon$ we update $i^n(\omega)$ by taking $i^n(\omega)-1$ ensuring that fractional distance exceeds the lower bound $\epsilon$ and changing $p^n(\omega)$ by $p^n(\omega)+1$. Consequently, 
\begin{equation}  \label{eq:bound_p_FT}
\epsilon < p^{n}(\omega) \leq 1+ \epsilon \,, \quad 0 \leq n \leq N, \ \omega\in{\Omega}.
\end{equation}

\begin{figure}[h]
\begin{center}
	\includegraphics[width=12cm]{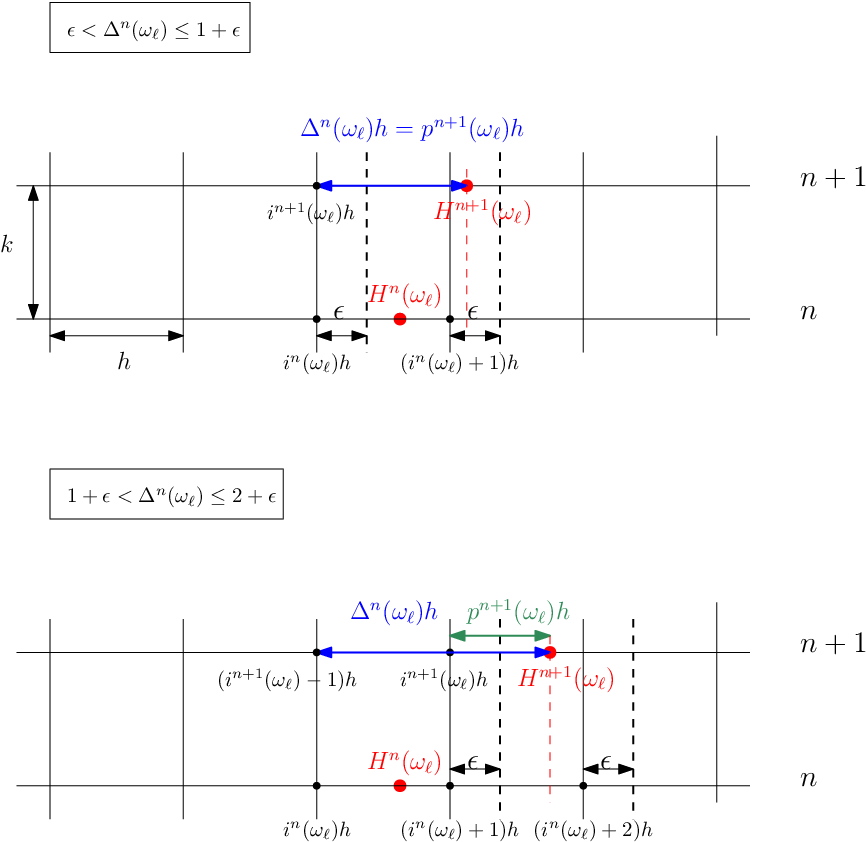}
\caption{The two possible scenarios where a sample of the front, $H^n(\omega_\ell)$, can be placed at time $t^n$. The sample fractional spatial distance, $\Delta^{n}(\omega_\ell)$, is computed by \eqref{eq:p_FT}. }
\label{fig:FTracking_cases}
\end{center}
\end{figure}

Taking into account \eqref{eq:p_FT} two scenarios depicted in Figure \ref{fig:FTracking_cases} can occur:
\begin{enumerate}
\item[a)]   $\epsilon <  \Delta^n(\omega) \leq 1+ \epsilon$. \\
In this case,  $i^{n}(\omega)\,h+ \epsilon < H^{n+1}(\omega) \leq (i^{n}(\omega)+1)\, h + \epsilon$, and then we take $i^{n+1}(\omega) = i^n(\omega)$. The fractional distance at level $n+1$ is $p^{n+1}(\omega)=\Delta^n(\omega)$.
\item[b)] $1+\epsilon<  \Delta^n(\omega) \leq 2+\epsilon$. \\
In this case , $(i^{n}(\omega)+1)\,h + \epsilon < H^{n+1}(\omega) \leq (i^{n}(\omega)+2)\, h+ \epsilon$, and then we take  $i^{n+1}(\omega) = i^n(\omega)+1$. The new fractional distance at level time $n+1$ is $p^{n+1}(\omega)=\Delta^n(\omega)-1$.

Note that here a new mesh spatial point is added, $r_{i^{n+1}(\omega)}=(i^{n}(\omega)+1)\, h$, and the numerical solution s.p. at this new point, $u^{n+1}_{i^{n+1}(\omega)}(\omega)$, has to be given. For this purpose we use a  Lagrange interpolation passing through the points $r_{i^{n}(\omega)-1}$, $r_{i^n(\omega)}$ and $H^{n+1}(\omega)$ if this quadratic interpolation results a positive value or the linear interpolation passing through the points $r_{i^n(\omega)}$ and $H^{n+1}(\omega)$ in other case, as follows
\begin{eqnarray}
u_{i^{n+1}(\omega)}^{n+1}(\omega) =
\left\{ \begin{array}{lll}
I_2(n+1;\omega), & \textrm{if}  & I_2(n+1;\omega)>0\,,\\
\\
I_1(n+1;\omega), &\textrm{if}  &  I_2(n+1;\omega)\leq 0\,,
\end{array} \right.
\end{eqnarray}
where
\begin{eqnarray}
I_2(n+1;\omega) & = & \frac{1-\Delta^{n}(\omega)}{1+\Delta^{n}(\omega)} u_{i^{n}(\omega)-1}^{n+1}(\omega) + 2 \frac{\Delta^{n}(\omega)-1}{\Delta^{n}(\omega)} u^{n+1}_{i^{n}(\omega)}(\omega)\,, \label{eq:I_2}\\
I_1(n+1;\omega) & = & \left(1-\frac{1}{\Delta^{n}(\omega)} \right) u_{i^{n}(\omega)}^{n+1}(\omega)\,. \label{eq:I_1}
\end{eqnarray}

\end{enumerate}

Other stability condition arrives from imposing the positivity of the coefficient $\tilde{B}_j^n(\omega)$   defined by \eqref{eq:CoeficientsAjnBjn_Cjn_FT}, to guarantee the positivity of the population solution s.p. in the interior points up to the last but one. It is easy to show that   $\tilde{B}_j^n(\omega) \geq 0$  under condition
\begin{equation}    \label{eq:bound2_k_FT}
k \leq \frac{h^2}{2d_2+ \left| \alpha_1-\beta_2\,P_0  \right| \,h^2}\,,
\end{equation}
where $d_2$ and  $P_0$ were defined in \eqref{eq:bound} and \eqref{eq:Def_elem_Theorem2}, respectively,  and the bounds $\alpha_1$  and $\beta_2$ were given in \eqref{eq:kappa_bound}.

Furthermore, the positivity of the numerical solution s.p. at the last interior point is guaranteed if the coefficient of $u^n_{i^n(\omega)}(\omega)$ in \eqref{eq:lastInteiorPoint}  is positive. Taking into account that $p^n(\omega) > \epsilon$, see \eqref{eq:bound_p_FT}, a straightforward calculus leads to the stability condition 
\begin{equation}  \label{eq:bound3_k_FT}
 k \leq \frac{\epsilon \ h^2 \ i^0}{d_2\, (2i^0+1-\epsilon)}\,,
\end{equation}
with $d_2$ defined in \eqref{eq:bound}.

Summarizing from conditions \eqref{eq:bound1_k_FT}, \eqref{eq:bound2_k_FT}--\eqref{eq:bound3_k_FT}, a stability condition for the RFDS-FT \eqref{eq:scheme_E1}--\eqref{eq:Front_FT}  has been found in order to guarantee the m.s. stability of the pairwise solution s.p.'s $\left\{ u_{j}^n(\omega),H(t^n;\omega) \right\}$:
\begin{equation}   \label{eq:StabilityCondition_FTmethod}
k \leq \min \left\{   \frac{h}{\eta_0\, \left| u^{\prime}_0(H_0) \right|}\,, \  \frac{h^2}{2d_2+ \left| \alpha_1-\beta_2\,P_0  \right| \,h^2}\,,\ \frac{\epsilon \ h^2 \ i^0}{d_2\, (2i^0+1-\epsilon)} \right\}\,.
\end{equation}

\subsection{Algorithm for computing the numerical  solution s.p.'s for FT method}   \label{subsec:Algorithm_FT}
Algorithm \ref{algo:FT} summarizes the procedure of the numerical solution s.p. employing the FT method for every realization $\omega_{\ell} \in \Omega$. 
\RestyleAlgo{ruled}
\begin{algorithm}
	\caption{Explicit RFDS-FT for the Stefan problem \eqref{eq:PDE}--\eqref{eq:BC} for a fixed $\omega_\ell\in \Omega$}\label{algo:FT}		
	\KwData{$D(\omega_\ell)$, $\eta(\omega_\ell)$, $\alpha(r)$, $\beta(r)$, $T$, $H_0$, $u_0$, $M$, $N$, $0 < \epsilon <1$}
	\KwResult{$r$, $u(r,T;\omega_\ell)$, $H(T;\omega_\ell)$}
	$h \gets H_0/M$\\
	$k \gets T/N$ \text{verifying the stability condition \eqref{eq:StabilityCondition_FTmethod}}\\
	$H^0 \gets H_0$\\
	$p^0\gets 1$\\
	$i^0(\omega_\ell) \gets M-1$ \\
	\For{$j = 0, \ldots, M-1$}{
		$r_j \gets jh$\\
		$u^0_j = u_0(r_j)$
	}
	\For{$n = 0, \ldots, N-1$}{
		\For{ $j =1, \ldots, i^{n}(\omega_\ell)-1$}{
			\text{Compute coefficients $\tilde{A}_j^n(\omega_\ell)$, $\tilde{B}_j^n(\omega_\ell)$ and $\tilde{C}_j^n(\omega_\ell)$ by \eqref{eq:CoeficientsAjnBjn_Cjn_FT}}\\
$u_j^{n+1}(\omega_\ell) \gets  \tilde{A}_j^n(\omega_\ell)\, u_{j-1}^n(\omega_\ell) + \tilde{B}_j^n(\omega_\ell)\, u_j^n(\omega_\ell) + \tilde{C}_j^n(\omega_\ell) \,u_{j+1}^n(\omega_\ell)$ \text{by \eqref{eq:scheme_E1}}
		}
		$u_0^{n+1}(\omega_\ell) \gets (4u_1^{n+1}(\omega_\ell) - u_2^{n+1}(\omega_\ell)) / 3$ \text{by} \eqref{eq:u0_FT}\\
		Compute $u_{i^n(\omega_\ell)}^{n+1}(\omega_\ell)$  by \eqref{eq:lastInteiorPoint}\\
$\Delta^n(\omega_\ell)  \gets p^n(\omega_\ell) + \frac{k \, \eta(\omega_\ell)}{h^2} \left(\frac{p^n(\omega_\ell)+1}{p^{n}(\omega_\ell)} \,u_{i^n(\omega_\ell)}^n(\omega_\ell)- \frac{p^n(\omega_\ell)}{p^{n}(\omega_\ell)+1}\,u_{i^n(\omega_\ell)-1}^n(\omega_\ell) \right)$
\text{by \eqref{eq:p_FT}}\\
		$H^{n+1}(\omega_\ell) \gets (\Delta^n(\omega_\ell)+i^n(\omega_\ell))\,h$\\
		\If{$\Delta^{n}(\omega_\ell) > 1 + \epsilon$}{ \text{Add new point,} $r_{i^{n+1}(\omega_\ell)}$: \\
		$r_{i^{n+1}(\omega_\ell)} \gets r_{i^n(\omega_\ell)}+ h$\\
		$i^{n+1}(\omega_\ell) \gets i^n(\omega_\ell)+1$\\
		\text{Approximate the unknown $u_{i^{n+1}(\omega_\ell)}^{n+1}(\omega_\ell)$ computing by a Lagrange interpolation passing through the points} $r_{i^{n}(\omega_\ell)-1}$, $r_{i^n(\omega_\ell)}$ \text{and} $H^{n+1}(\omega_\ell)$:	\\
		$u_{i^{n+1}(\omega_\ell)}^{n+1}(\omega_\ell) \gets I_2(n+1;\omega_\ell)$ \text{using} \eqref{eq:I_2}\\
	\If{$u_{i^{n+1}(\omega_\ell)}^{n+1}(\omega_\ell) <0$}{\text{use a linear interpolation to recalculate the value}\\
	$u_{i^{n+1}(\omega_\ell)}^{n+1}(\omega_\ell) \gets I_1(n+1;\omega_\ell)$ \text{using} \eqref{eq:I_1} \\
		}
			$p^{n+1}(\omega_\ell) \gets \Delta^{n}(\omega_\ell)-1$ 
		}
$p^{n+1}(\omega_\ell)=\Delta^n(\omega_\ell)$	\\		}
	\end{algorithm}


\section{Monte Carlo technique}  \label{sec:MCmethod}

Iterative methods such as finite difference schemes become unsuitable in the random scenario because of their unworkable complexity.  It is caused by the storage accumulation of the symbolic computation of the involved s.p. solutions at the intermediate levels, see \citep{CasabanCompanyJodar2021}. To avoid these drawbacks we use a combination of the sample m.s. approach with Monte Carlo method, see \cite{Kroese2011}. This method allows us to compute efficiently the statistical moments of the approximating solution stochastic process and the stochastic moving boundary solution.

In this way we take a given number $K$ of realizations for the random parameters involve in the RFDS-(FF/FT) and compute the solutions $\left\{ \left(u^{n}(\omega_\ell),H^{n}(\omega_\ell) \right), \ 1 \leq \ell \leq K \right\}$ using Algorithm \ref{algo:EFDM} and Algorithm \ref{algo:FT}. 
The mean and the standard deviation of the solutions are computed simultaneously by using the following well-known formulae 
\begin{equation}
    \mu[u] = \frac{1}{K} \sum_{\ell = 1}^K u_{\ell}, \quad \sigma[u] = \sqrt{\frac{1}{K}\sum_{\ell = 1}^K\left( u_{\ell} \right)^2- \left(\mu[u]\right)^2},
\end{equation}
where $u_{\ell}$, $1\leq \ell \leq K$, denotes the solution for a sample realization $\omega_{\ell}$. Similar formulae are used for the random free-boundary. Algorithm \ref{algo:MC} summarizes the procedure of the computation of the statistical moments of the numerical  solution s.p. for the RPDE problem \eqref{eq:PDE}--\eqref{eq:BC}. Note that the same algorithm is used for computation of the statistical moments of the random free-boundary.

\RestyleAlgo{ruled}
\begin{algorithm}
	\caption{Monte Carlo technique for the statistical moments of the numerical solution s.p. for RPDE problem \eqref{eq:PDE}--\eqref{eq:BC}}\label{algo:MC}
	\KwData{$K$}
	\KwResult{$\mu[u]$, $\sigma[u]$}
	$\text{SUM} \gets 0$ \;
	$\text{SUM2} \gets 0$\;
	\For{ \text{each realization } $\omega_{\ell}$, $1\leq \ell \leq K$}{
		\text{Set r.v. parameters for the problem \eqref{eq:PDE}--\eqref{eq:BC}}\\
		$u_{\ell}  \gets \text{numerical solution of \eqref{eq:PDE}--\eqref{eq:BC} for } \omega_{\ell} \text{ by Algorithm \ref{algo:EFDM} (FF) or Algorithm  \ref{algo:FT}} $  (FT)\\
		$ \text{SUM} \gets \text{SUM} + u_{\ell} $\\
		$\text{SUM2} \gets  \text{SUM2} + \left( u_{\ell} \right)^2$
	}
	$\mu[u] \gets \text{SUM} / K$\\
	$\sigma[u] \gets  \sqrt{\text{SUM2} / K - \left(\mu[u]\right)^2}$
\end{algorithm}

As a remark for the implementation of the code, note that in the case of the random FF method the number of the grid points is fixed for all the time iterations and uniform for all the sample realizations of the r.v.'s $D(\omega)$ and $\eta(\omega)$.  However, in the case of the random FT method  the final count of interior grid points depends on the r.v.'s $D(\omega)$ and $\eta(\omega)$, thus resulting in variability across realizations. Therefore, to facilitate the computation of the moments at the final time $T$, we standardize by $I(\ell)$ the maximum number of the interior points up to $T$ for each realization $\omega_{\ell}$. As $I (\ell)$ is not known beforehand, we designate $I_{\textrm{max}}$ to be large enough such that it satisfies $I_{\textrm{max}} \geq \max_{1\leq \ell\leq K} I(\ell)$. Consequently,

\begin{equation}
	u_j(\omega_\ell) = \begin{cases}
	u_j(\omega_\ell), & 0 \leq j\leq I(\ell),\\
	0, & I(\ell) \leq j \leq I_{\textrm{max}},
\end{cases}
\end{equation}

To speed up the computational process, which becomes particularly crucial when dealing with large-scale Monte Carlo simulations, we have employed a parallel computing environment. In this setup, each processor independently computes the solutions for every realization of the random variables. The simultaneous execution of computations across different processors significantly boosts the efficiency and speed of Monte Carlo simulations.

For parallel computing, MATLAB provides a specialized tool known as the Parallel Computing Toolbox. This toolbox includes an improved variant of the \texttt{for} loop, referred to as \texttt{parfor}. This function enables iterations of a particular code block to be executed a specified number of times in parallel, thus offering a practical and efficient way to leverage the advantages of parallel computing in Monte Carlo simulations.


\section{Numerical examples}  \label{sec:NumericalExamples}
In this section, we compare the numerical results from FF and FT methods  as defined in Algorithms \ref{algo:EFDM} and \ref{algo:FT}. First, we use these methods on one specific realization of the random free-boundary diffusive logistic model \eqref{eq:PDE}--\eqref{eq:ConditionsCI}, treating it like a predictable, or deterministic, case. Next, we compare results from the random versions of the FF and FT methods when used on the same random model \eqref{eq:PDE}--\eqref{eq:ConditionsCI}.

Lastly, we analyse the numerical convergence of  the random FF method developed by means of the study of the convergence of both the Monte Carlo method and the RFDS given by \eqref{eq:scheme_coef}--\eqref{eq:scheme_BC}.

\subsection{Comparison between the random FF and FT methods}  
\label{subsection:comparisonFF-FT}

\begin{example}
	Deterministic case. 
\end{example}
Firstly, we compare both methods for a deterministic, constant scenario, using the parameters from Table \ref{table:Comparison_parameters} and setting both $D$ and $\eta$ to 1. For $M=50$ spatial points, the spatial step-sizes for the deterministic FF and FT methods become $h_{FF}=\frac{1}{M}=\frac{1}{50}$ and $h_{FT}=\frac{H_0}{M}=\frac{3}{50}$, respectively. Considering that the stability condition \eqref{eq:StabilityCondition_FTmethod} for the FT method ($k<8.9543\mathrm{e}-04$) is stricter than the one for the FF method ($k<0.0016$), we adopt a time step-size of $k=8.0\mathrm{e}-04$.

\begin{figure}[h]
	\begin{center}
	
			\includegraphics[width=\textwidth]{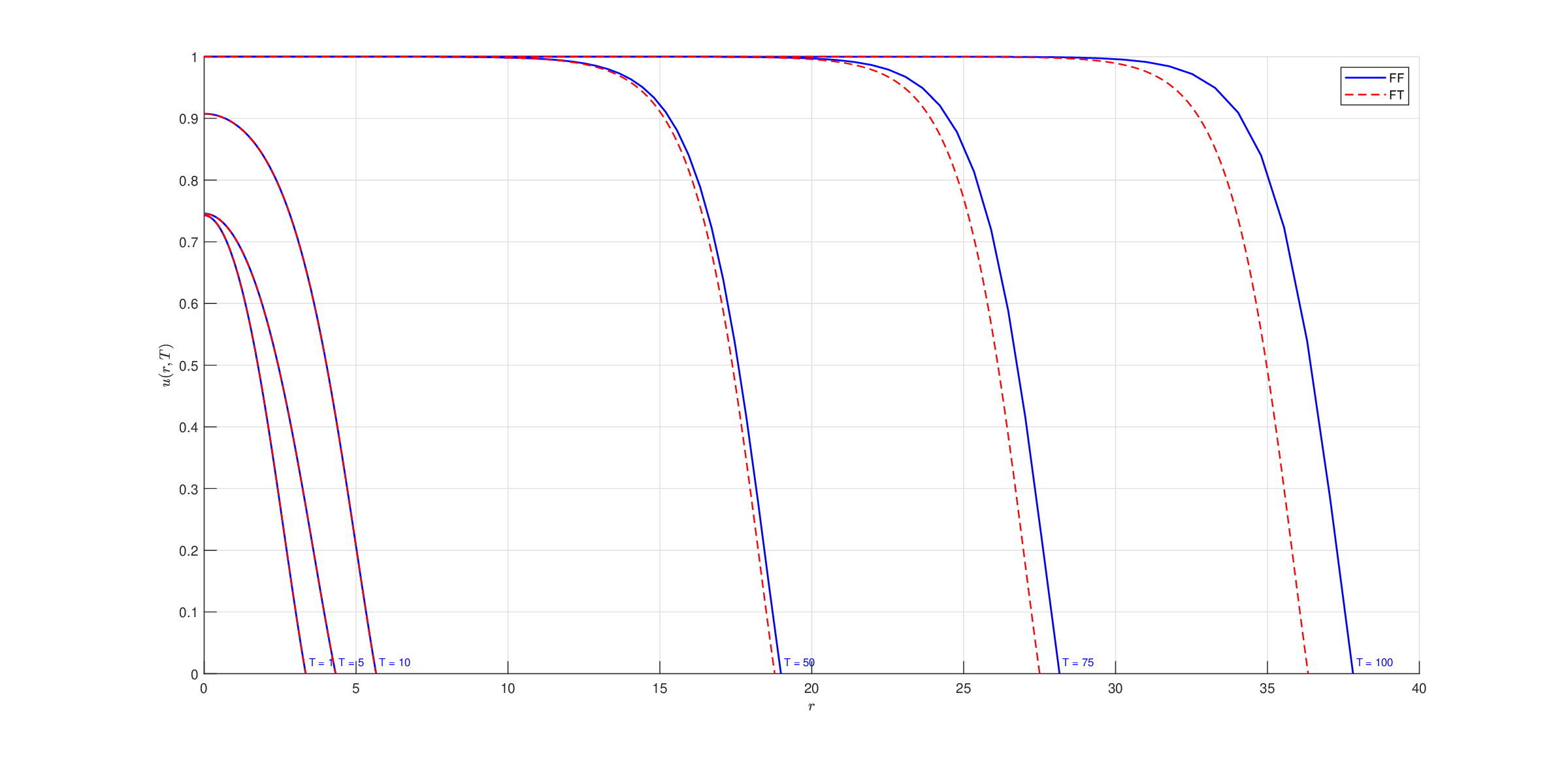}
	
		\caption{Numerical solution for the deterministic case of the problem \eqref{eq:PDE}--\eqref{eq:ConditionsCI}   derived from the FF method (blue solid line) and the FT method (red dashed line) across various values of $T=\{1, 5,10,50,75,100 \}$.}
		\label{fig:determ_T}
	\end{center}
\end{figure}

Figure \ref{fig:determ_T} presents the numerical solution derived from the proposed methods across various values of $T$. The discrepancies between the two methods, particularly evident at $T>50$, are primarily due to the FF method's loss of accuracy as the domain size increases.  The FF method retains the number of grid points and the step-size within the transformed fixed spatial domain $[0,\, 1]$. However, when the inverse transform is applied to revert to the original variables, the step-size $\Delta r$ escalates with an increasing domain size $[0, \, H(T)]$, leading to a decline in the method's accuracy. In contrast, the FT method's ability to maintain the step-size and add new grid points as the domain expands ensures its consistency in terms of accuracy. Therefore, we can infer that the FF method is more effective for smaller $T$ values, whereas the FT method is the preferable choice for larger $T$ values.

\begin{figure}[h]
	\begin{center}
		\begin{minipage}[t]{0.45\textwidth}
			\includegraphics[width=\textwidth]{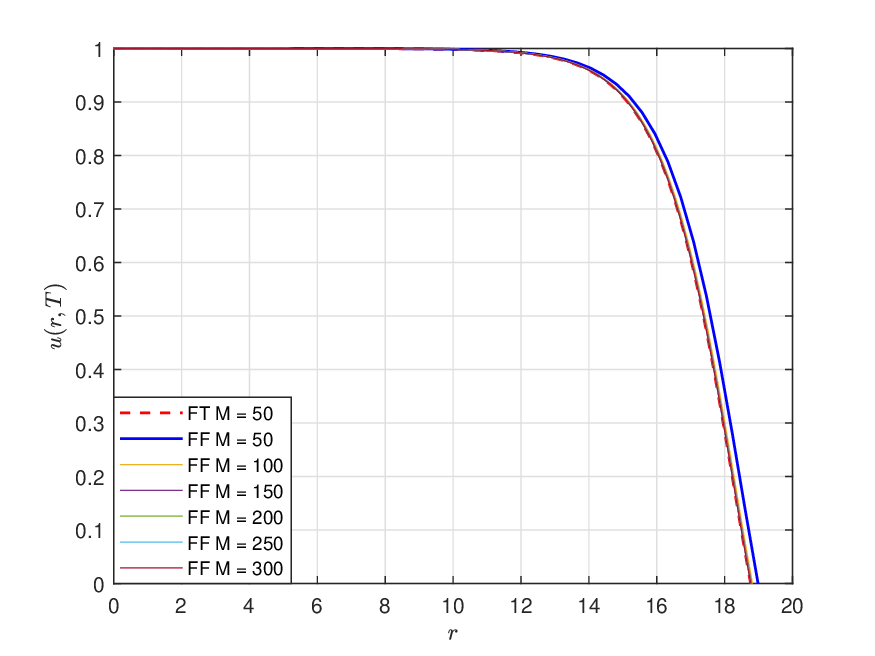}
		\end{minipage}
		\begin{minipage}[t]{0.45\textwidth}
			\includegraphics[width=\textwidth]{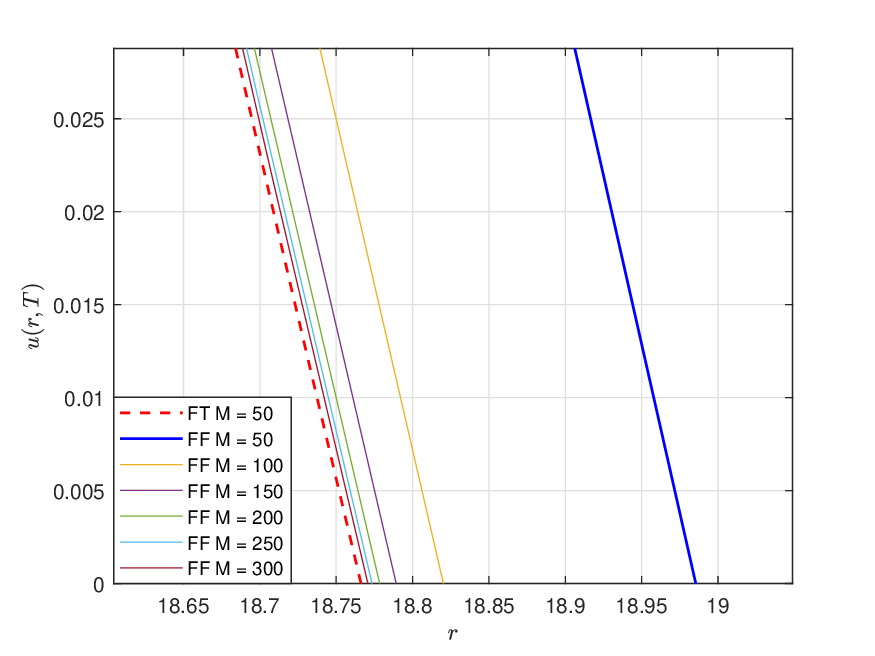}
		\end{minipage}
		\caption{Numerical solution for the deterministic version of problem \eqref{eq:PDE}--\eqref{eq:ConditionsCI}, comparing the FF method with different numbers of grid points $M$, and the FT method using $M=50$ grid points.The plot on the right shows a zoomed-in view for $r$ close to $H(T)$ at $T=50$.}
			 \label{fig:determ_M}
	\end{center}
\end{figure}

To illustrate this observation, we fix $T = 50$ (a point at which the discrepancy becomes noticeable) and modify the number of grid points $M$ in the FF method. The results are visualized in Figure \ref{fig:determ_M}. In the plot on the right, we provide a zoomed-in view of the solutions close to $r = H(T)$. This helps highlight the convergence to the results from the FT method. The total number of nodes at the final moment $t=T$ for the FT method is $i^{N}+2 = 315$,  where $Nk=T=50$. Table \ref{table:CPUs_FFMs_FT} collets the CPU,s as well as the real time lapsed corresponding to these CPU,s for the computations carry out in Figure \ref{fig:determ_M}. Computations have been carried out by Matlab$^{\copyright}$ software version R2019b Update 3 for Windows 10Pro (64-bit) AMD Ryzen Threadripper 2990WX 32-Core Processor, 3.00 GHz. The timings have been computed using cputime function of Matlab$^{\copyright}$ (CPU time spent).

\begin{table}[!h]
\begin{center}
\begin{tabular}{@{}cccccc}  \hline
   &  FF  & FF & FT & FT 
\\
$M$ &  CPU-time (s)  & real time (s)  &  CPU-time (s)  & real time (s) 
\\
\hline 
  $50$ &   $0.4063$       &   $0.27$   & $0.8906$   & $0.50$
 \\
 $100$     &   $1.2500$   & $1.13$  &  - & -
 \\
 $150$   &   $3.0625$    & $2.76$    & - & -
\\      
 $200$   &   $5.4688$    & $5.23$    & - & -
 \\
  $250$   &   $9.6406$    & $9.32$   & - & -
  \\
$300$   &   $16.0313$    &  $15.48$   & - & -
\\
\hline
\end{tabular}
\caption{CPU time seconds and their corresponding real time in seconds spent to compute the numerical results depicted in Figure \ref{fig:determ_M}. }\label{table:CPUs_FFMs_FT}
\end{center}
\end{table}

With our findings from the deterministic case in mind, we now move on to studying the random case, which is the main focus of our current research.

\begin{example}
	Random case.
\end{example}
This example is a comparative study of the deviations of the numerical solution s.p.'s, $\{ (u(r,t;\omega),H(t;\omega))\}$, between the random FF and FT methods. The values of the parameters considered in this study for the model \eqref{eq:PDE}--\eqref{eq:ConditionsCI}  are collected in Table \ref{table:Comparison_parameters}. Note that we have chosen  truncated distributions for the s.v.'s $D(\omega)$ and $\eta(\omega)$ because  of the boundedness conditions \eqref{eq:bound} and \eqref{eq:bound_eta}. We distinguish between the constant case and the variable case for the parameters $\alpha$ and $\beta$ involve in the model. 
 
In order to guarantee the stability of the solutions s.p. generated for both methods and in the spirit of comparing them, we take a common stability condition. Taking $M=50$ spatial points, the spatial step-sizes for the random FF and FT methods are, $h_{FF}=\frac{1}{M}=\frac{1}{50}$, and $h_{FT}=\frac{H_0}{M}=\frac{3}{50}$, respectively.  Due to the fact that the stability condition \eqref{eq:StabilityCondition_FTmethod} of random FT method ($k<7.4619\mathrm{e}-04$) is more restrictive than one imposed on the random FF method ($k<0.0013$), we consider the following time step-size $k=7.0\mathrm{e}-04$.

\begin{table}[h]
    \centering
    \begin{tabular}{ll}
    \hline
    Constant case & Variable case\\
    \hline
         $D(\omega) \sim \mathcal{N}_{[0.8,\, 1.2]}(1, \, 0.1)$ &   $D(\omega) \sim \mathcal{N}_{[0.8,\, 1.2]}(1, \, 0.1)$ \\
         \medskip
         $\eta(\omega) \sim \mathcal{B}e_{[1.6,\, 2.4]}(2, 4)$ & $\eta(\omega) \sim \mathcal{B}e_{[1.6,\, 2.4]}(2, 4)$\\
         \medskip
         $\alpha=1$ & $\alpha(r)=\frac{2r+3}{2r+2}$\\
         \medskip
         $\beta=1$ & $\beta(r)=\frac{2r+1}{2r+2}$ \\
         \medskip
         $H_0=3$ & $H_0=3$\\
         $u_0(r) = \cos \left( \frac{\pi\,r}{6}\right)$, \ {\small $0 \leq r \leq H_0$} & $u_0(r) = 1 - \left(\frac{r}{H_0}\right)^2$, \ {\small $0 \leq r \leq H_0$}\\
         \hline
    \end{tabular}
    \caption{Data used in the comparison between the random FF and FT methods for the constant and variable cases of the parameters $\alpha$ and $\beta$.}
    \label{table:Comparison_parameters}
\end{table}

Table \ref{table:Comparison_CPU_variable} reports the CPU-time seconds and the corresponding elapsed real time for the variable case.  As we can observe, both methods provide competitive times and slightly lower for the random FF method.
 
\begin{table}[!h]
\begin{center}
\begin{tabular}{@{}cccccc}  \hline
 &   & & $[\mu$/$\sigma]$ $\left(u_K,H_K\right)$ \ \ { \small for variable case: $\alpha(r)$, $\beta(r)$}  & \\
 \hline
  &  FF  & FF & FT & FT 
\\
$K$ &  CPU-time (s) & real time (s)  &   CPU-time (s) & real time (s) 
\\
\hline 
  $25$ &   $0.5625$ &$0.33$   &  $0.5781$  & $0.36$
 \\
 $50$     &   $0.4375$    & $0.43$  &  $0.8125$ &  $0.5$
 \\
 $100$   &   $0.6563$    &  $0.48$  & $1.0469$ & $0.56$
\\      \hline
\end{tabular}
\caption{CPU time seconds and their corresponding real time in seconds spent to compute the $K$-Monte Carlo simulations for the mean, $\mu$, and the standard deviation, $\sigma$, of solutions $u_K(r_j,T;\omega_\ell)$ and $H_k(T;\omega_\ell)$, $1 \leq \ell \leq K$, at  $T=1$ for the variable case of parameters of Table \ref{table:Comparison_parameters}. }\label{table:Comparison_CPU_variable}
\end{center}
\end{table}
We have computed the maximum of the relative deviations between the solutions of the FF method,  $u^{FF}(r_j,T;\omega_\ell)$,  and the FT method, $u^{FT}(r_j,T;\omega_\ell) $, for all the realizations $\omega_\ell$,  $1 \leq \ell \leq K$,  in the set $K=\{25,50,100 \}$ by means

\begin{equation}   \label{MaxRErr-FFvsFT_u}
\textrm{RelErr} \left( u_{K}^{FF},u_{K}^{FT} \right) =  \max_{ 0  \leq j \leq M, \, 1 \leq \ell \leq K } 
\left| \frac{u^{FF}(r_j,T;\omega_\ell) - u^{FT}(r_j,T;\omega_\ell) }{u^{FF}(r_j,T;\omega_\ell)} \right|   
\end{equation}

\begin{figure}[h]
\begin{center}
\begin{minipage}[t]{0.45\textwidth}
	\includegraphics[width=\textwidth]{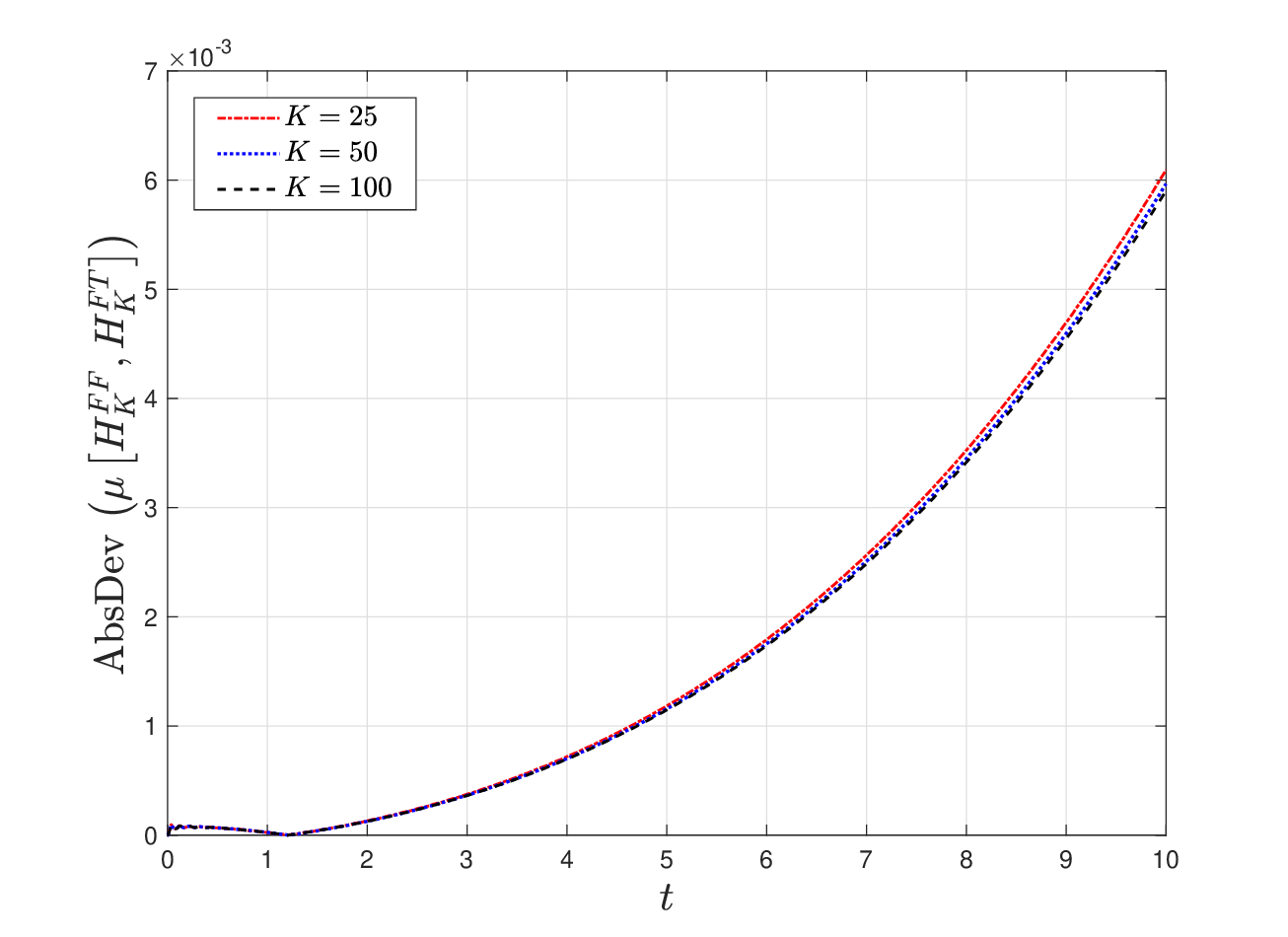}
\end{minipage}
\begin{minipage}[t]{0.45\textwidth}
	\includegraphics[width=\textwidth]{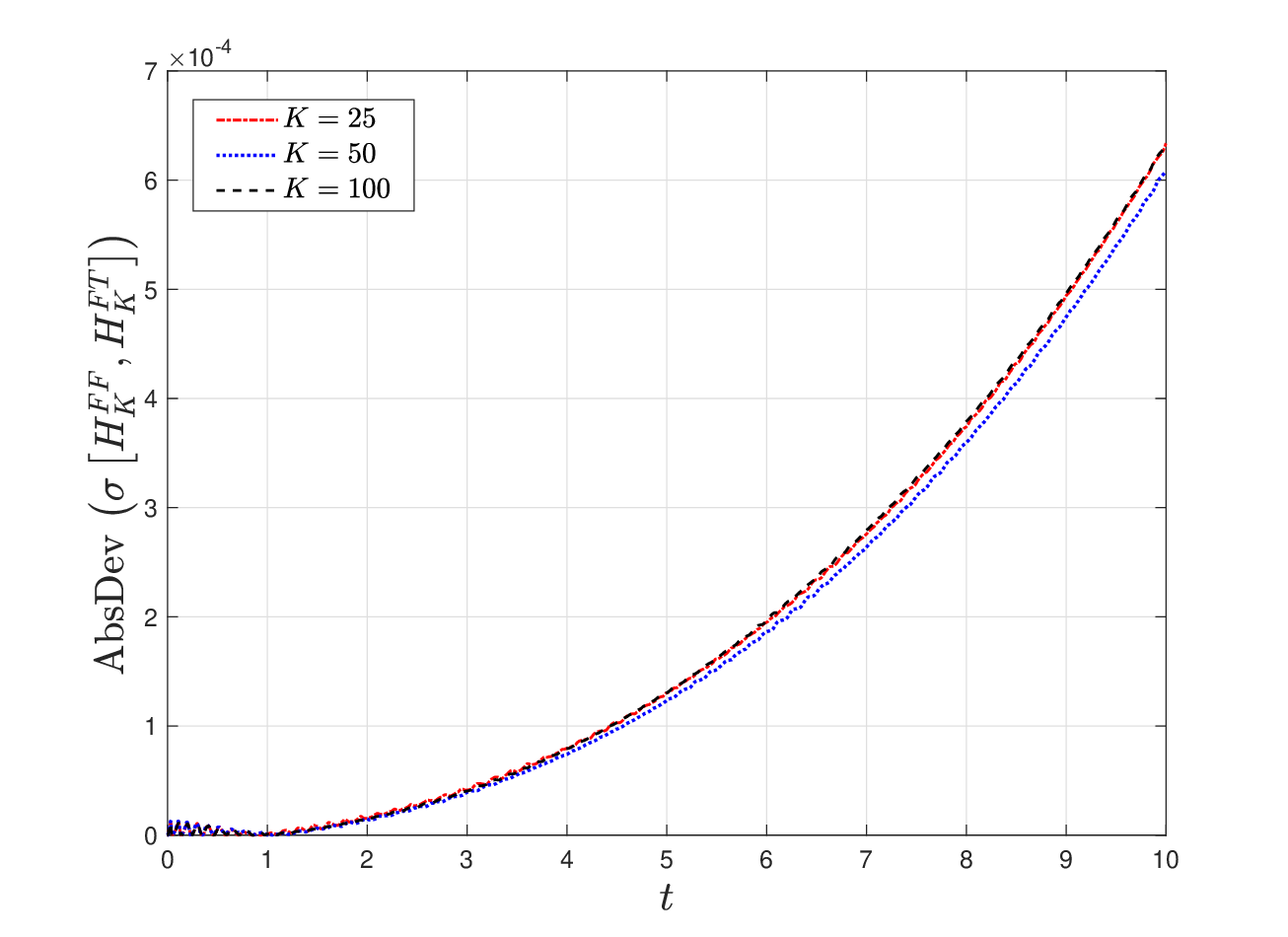}
\end{minipage}
\caption{Absolute deviations of the mean (left) and the standard deviation (right) of the random free boundary $H(T;\omega_{\ell})$, $1 \leq \ell \leq K$, computed  between random FF and FT method up to time $T=10$ by \eqref{AbsErr-FFvsFT_H} when parameters $\alpha$ and $\beta$ are constant (see Table \ref{table:Comparison_parameters}) and the sample $K$ realizations varies.  }
\label{fig:comparisonFFvsFT_constantCase}
\end{center}
\end{figure}

\begin{figure}[h]
\begin{center}
\begin{minipage}[t]{0.45\textwidth}
	\includegraphics[width=\textwidth]{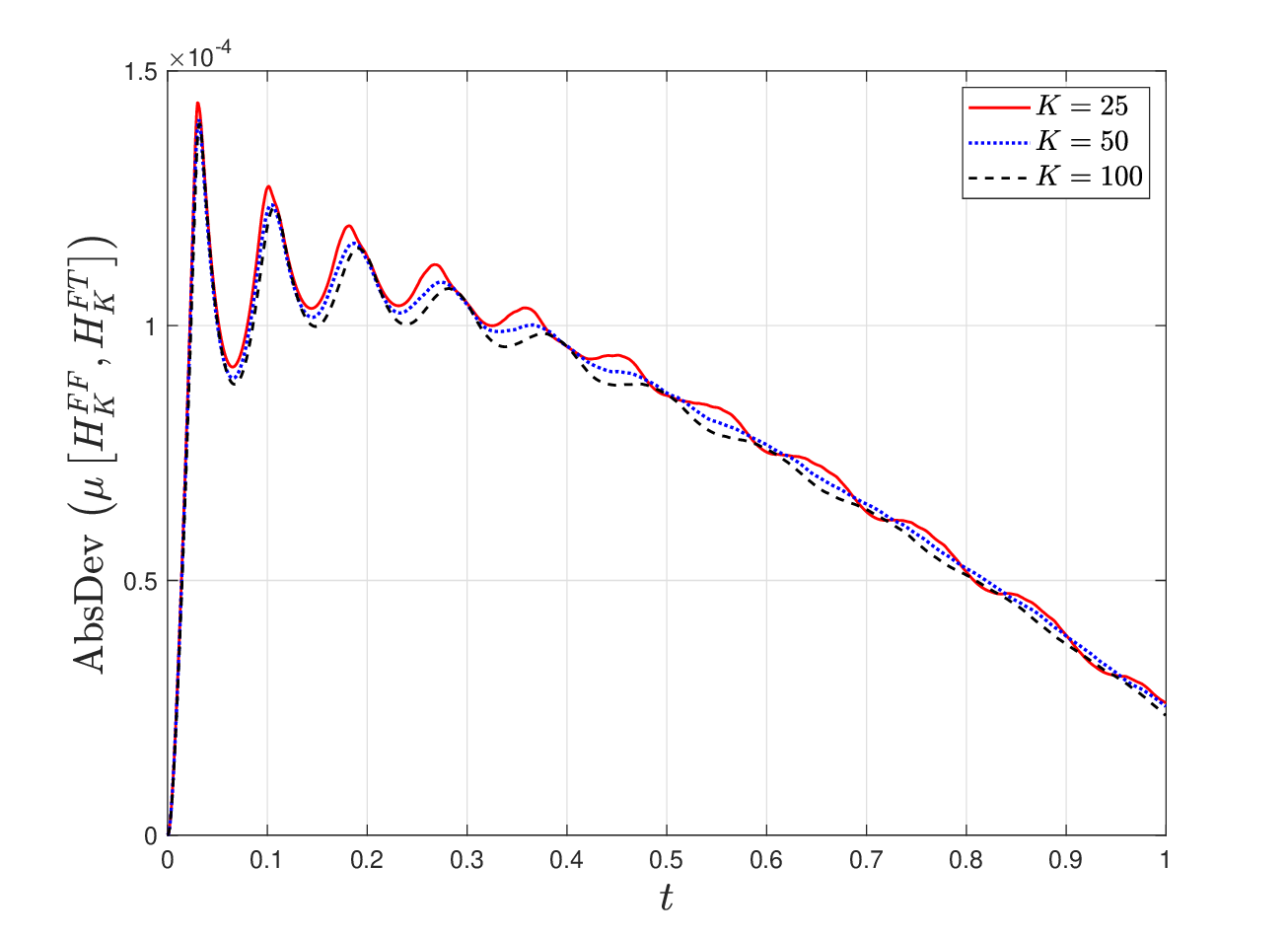}
\end{minipage}
\begin{minipage}[t]{0.45\textwidth}
	\includegraphics[width=\textwidth]{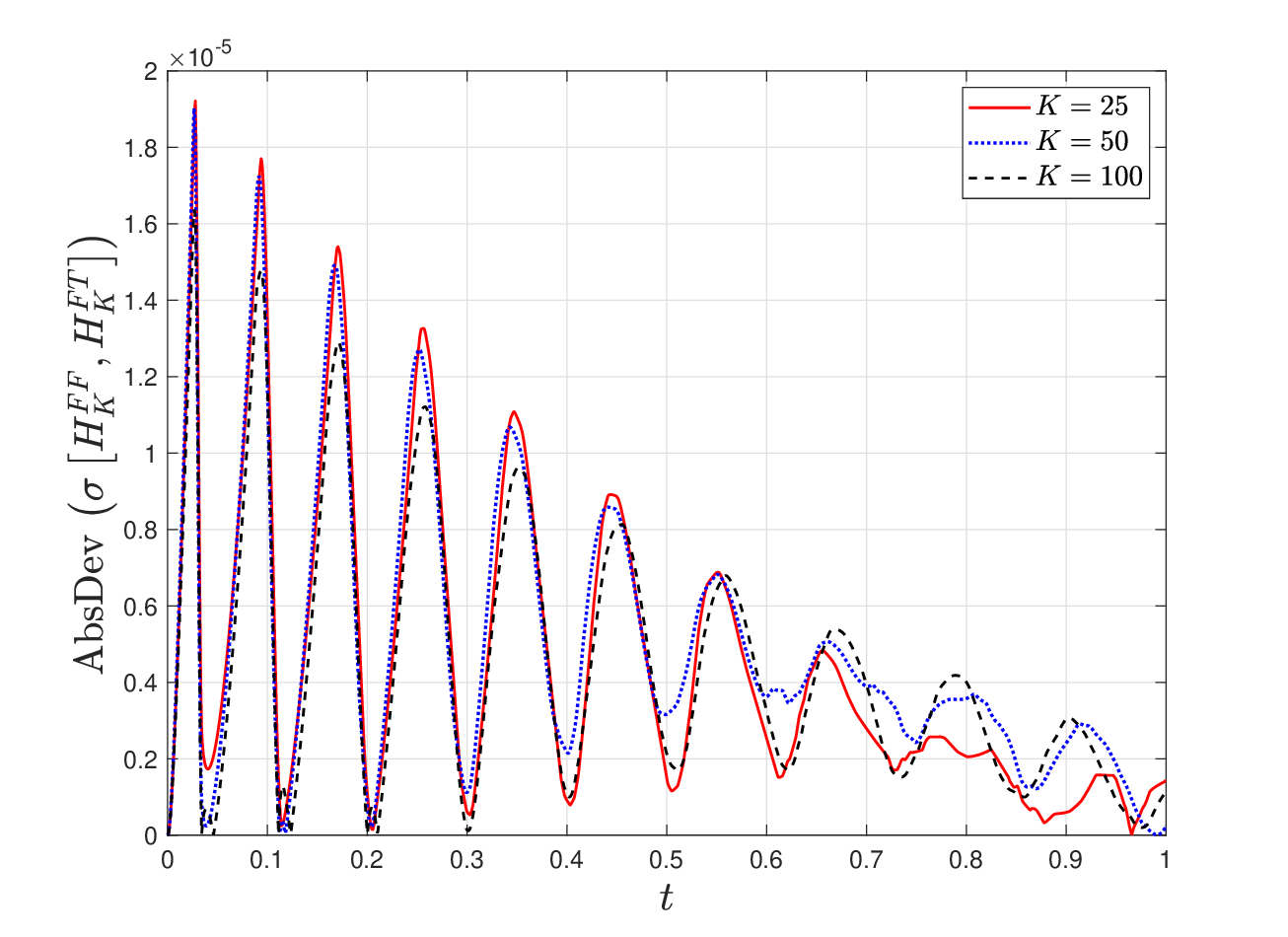}
\end{minipage}
\caption{Absolute deviations of the mean (left) and the standard deviation (right) of the random free boundary $H(T;\omega_{\ell})$, $1 \leq \ell \leq K$, computed  between random FF and FT method up to time $T=1$ by \eqref{AbsErr-FFvsFT_H} when parameters $\alpha$ and $\beta$ are variable (see Table \ref{table:Comparison_parameters}) and the sample $K$ realizations varies.}
\label{fig:comparisonFFvsFT_variableCase}
\end{center}
\end{figure}

In the comparison \eqref{MaxRErr-FFvsFT_u}, we take $M=50$ spatial points. But as the spatial grid, $r_j$, is not the same in both methods due to the nature of the methods,  firstly, we construct  the solution cubic spline curves generated by FT method, for each number of realizations of the set $K=\{25,50,100 \}$, and then we evaluate these solution curves in the spatial grid of the FF method generated up to $T$, that is, $r_j(\omega_\ell)=j\, \Delta r (\omega_\ell)=j\,\frac{H^{FF}(T;\omega_{\ell})}{M}$, $j=0,\ldots,M$. We denote $H^{FF}(T;\omega_{\ell})$  the localization of the moving front for the realization $\omega_{\ell}$ computed by FF method up to $T$ time. Table \ref{table:Comparison_ErrRel_u} collects the relative errors \eqref{MaxRErr-FFvsFT_u} for several simulations $K=\{25,50,100 \}$ at $T=10$ for the constant case of parameters of Table \ref{table:Comparison_parameters} and at $T=1$ time for the variable case. It is observed that both methods generate similar solutions regardless of the number of realizations $K$.
In the study of errors of the random free  boundary, $H(t;\omega)$, we use the following absolute deviations
\begin{eqnarray} \label{AbsErr-FFvsFT_H}
\begin{array}{lcl}
\textrm{AbsDev}  \left(\mu \left[ H^{FF}_{K},H^{FT}_{K} \right] \right) & = &  
\left| \, \mu_K \left[ H^{FF}(t;\omega_{\ell}) \right]  -  \mu_K \left[ H^{FT}(t;\omega_{\ell}) \right] \, \right|\,, \quad 0 \leq t \leq T, \ K  \ \text{fixed}, \\
\\
\textrm{AbsDev}  \left(\sigma \left[ H^{FF}_{K},H^{FT}_{K} \right] \right) & = &
  \left| \, \sigma_K \left[ H^{FF}(t;\omega_{\ell}) \right]  -  \sigma_K \left[ H^{FT}(t;\omega_{\ell}) \right] \, \right|\,, \quad 0 \leq t \leq T, \ K  \ \text{fixed} \,.
\end{array}
\end{eqnarray}  
Figure \ref{fig:comparisonFFvsFT_constantCase} and Figure \ref{fig:comparisonFFvsFT_variableCase} show that the approximations of both statistical moments of the free boundary  $\mu_K[H(t;\omega_{\ell})]$ and $\sigma_K[H(t;\omega_{\ell})]$, where $\mu_{K}$ and $\sigma_{K}$ denote the mean and the standard deviation respectively  for a number $K$ of realizations $\omega_{\ell}$, $1\leq \ell \leq K$, computed by means of the random FF and FT methods are close independently of the number of Monte Carlo simulations $K$. Note that in these figures we illustrate the full evolution of the absolute deviations of the statistical moments of the free boundary up to the times $T=1$ and $T=10$, respectively.

\begin{table}[!h]
\begin{center}
\begin{tabular}{@{}cccc}  \hline
$K$ &   $\textrm{RelErr} \left( u_{K}^{FF},u_{K}^{FT}\right)$ &   $\textrm{RelErr} \left( u_{K}^{FF},u_{K}^{FT} \right)$           \\
       &                {\small (constant case): $(\alpha,\beta)$)}   &         {\small (variable case: $(\alpha(r),\beta(r))$)}       
\\
 \hline
  $25$ &   $5.9649 \textrm{e}{-04}$   & $4.5482  \textrm{e}{-05}$
 \\
 $50$     &   $5.9900 \textrm{e}{-04}$   & $ 4.5220\textrm{e}{-05}$
 \\
 $100$   &   $6.3451 \textrm{e}{-04}$    & $ 4.6164 \textrm{e}{-05}$
\\      \hline
\end{tabular}
\caption{Relative errors, \eqref{MaxRErr-FFvsFT_u}, between the solutions of the FF method,  $u^{FF}$,  and the FT method, $u^{FT}$, at  time $T=10$ for the constant and variable case of parameters $\alpha$ and
 $\beta$ (see Table \ref{table:Comparison_parameters}), respectively, among number of MC realizations $K$. }\label{table:Comparison_ErrRel_u}
\end{center}
\end{table}

\begin{figure}[h]
	\begin{center}
					\includegraphics[width=\textwidth]{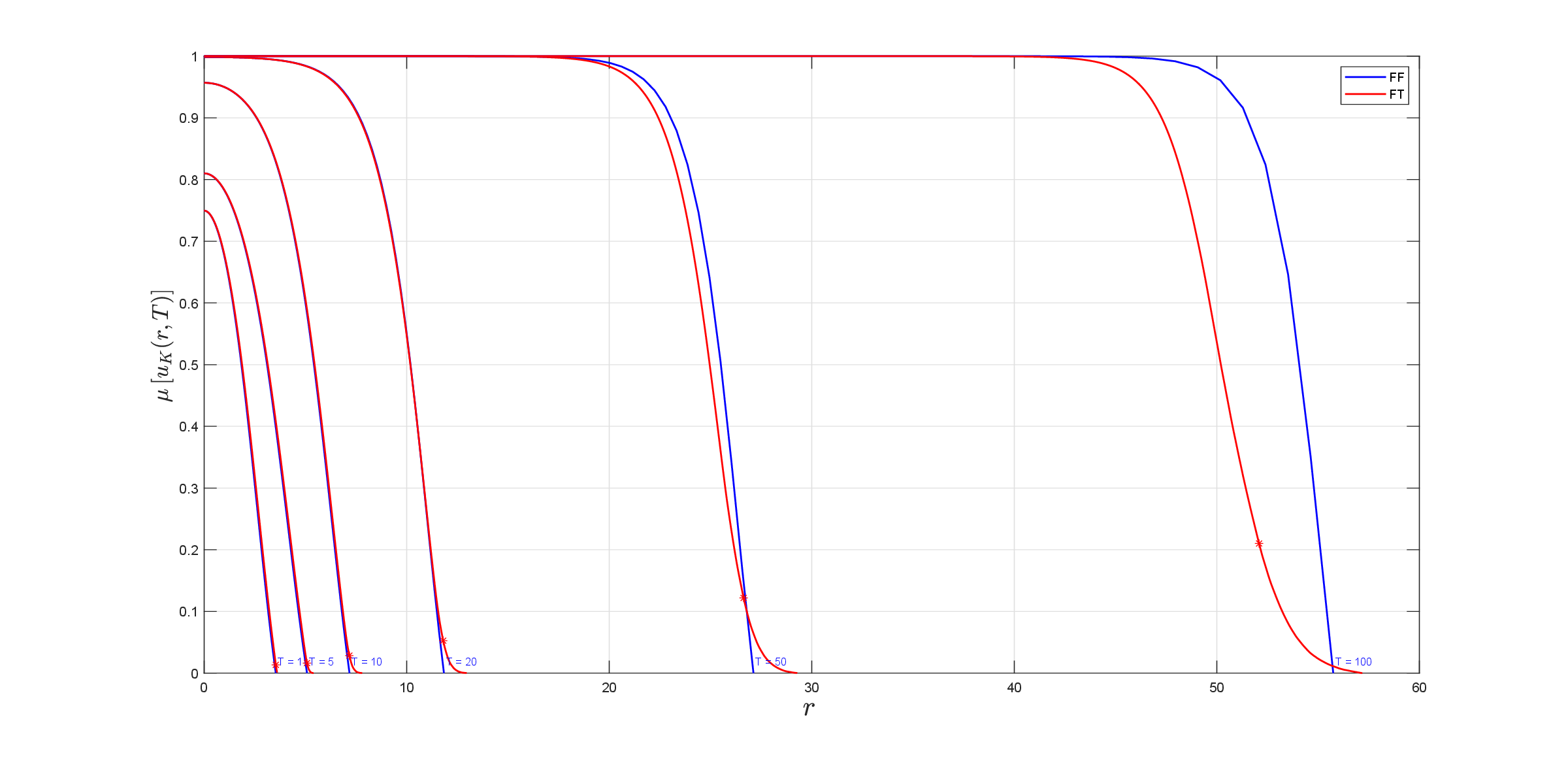}
				\caption{Mean of numerical solution s.p.'s  computed by the FF and FT methods over different time times $T=\{1,5,10,20,50,100 \}$. }
		\label{fig:random_u_T}
	\end{center}
\end{figure}

\begin{figure}[h]
	\begin{center}
		\includegraphics[width=\textwidth]{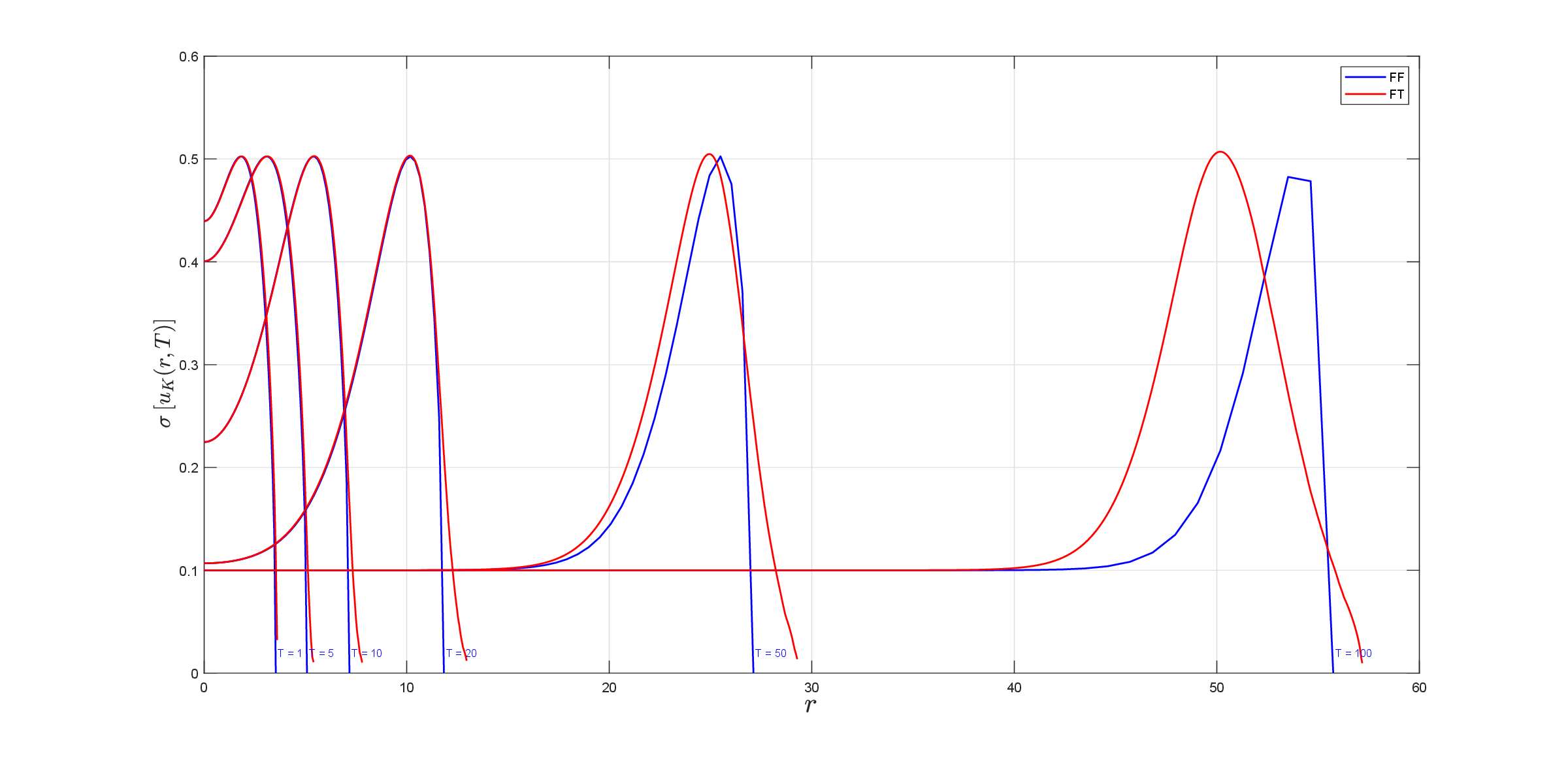}
		\caption{Standard deviation of numerical solution s.p.'s  computed by the FF and FT methods over different times $T=\{1,5,10,20,50,100 \}$. }
		\label{fig:random_u_T_sig}
	\end{center}
\end{figure}

Similar to the approach for the deterministic case, we compare the proposed method over different time spans $T$ for the constant case. We set $M = 50$ and $K=100$. The mean value and standard deviation of the numerical solution s.p.'s for various time intervals $T$ are presented in Figures \ref{fig:random_u_T} and \ref{fig:random_u_T_sig}, respectively. In Figure \ref{fig:random_u_T}, the red asterisks mark the positions of the average of $H(T;\omega_{\ell})$ for the corresponding $T$ calculated by the FT method. 
This point represents an inflection point on the curve of the mean. Beyond this value, the average accumulates non-zero values from the realizations, due to the varying positions of the free boundary in each realization. This variability leads to a convex shape in the curve. The average value drops to zero at the point where the free boundary $H(T;\omega_{\ell})$ reaches its maximum across all realizations. 
On the other hand, as a result of the front-fixing transformation \eqref{eq:landau}, the spatial variable becomes r.v.  for a fixed $t>0$, and one gets
\begin{equation}  \label{eq:rv_r}
r_j(t;\omega_\ell) = j \, \frac{H^{FF}(t;\omega_\ell)}{M}\,, \quad 0 \leq j \leq M\,.
\end{equation}
Thus the statistical moments of the approximate solution s.p. $u(r,t;\omega)$ are assigned to the mean of the spatial r.v. $\bar{r}_j(t)$, i.e.
\begin{equation}\label{eq:mu_r}
 \mu \left[ u^{FF} \left( \bar{r}_j(t^n),t^n;\omega_\ell \right) \right]=
 \mu \left[ v_j^n(\omega_\ell)  \right]\,, \quad \text{with} \ \ \ \bar{r_j}(t)=  \frac{j}{M} \, \mu \left[H^{FF}(t;\omega_\ell)\right]\,,
\end{equation}
where $\left\{v_j^n(\omega_\ell)\right\}$ is the sample numerical solutions of RFDS-FF \eqref{eq:scheme_coef}.
In particular, taking into account 
$v^{n}_M(\omega_\ell)=0$, $0 \leq n \leq N$, see \eqref{eq:scheme_BC}, the statistical moments assigned to the average to the moving front $H(t;\omega_{\ell})$ are identically zero
\[
 \mu \left[ u^{FF} \left( \bar{r}_M(t),t;\omega_\ell \right) \right]=0\,, \qquad   \sigma \left[ u^{FF} \left( \bar{r}_M(t),t;\omega_\ell \right) \right]=0\,.\]

 Formulae \eqref{eq:mu_r} also explain  the FF method's loss of accuracy as the domain size increases at each fixed realization. As in the deterministic case, this drawback can be avoided by taking a smaller spatial step-size for FF method.

As seen in Figure \ref{fig:random_u_T_sig}, the standard deviation of the numerical solutions is close to zero when the population approximates  the carrying capacity since only small variations across the different realizations have been found as expected. As we get closer to the free boundary, the standard deviation for the FF method also moves towards zero, as for all realizations the population density becomes zero when the spatial variable reaches the moving front. However, in contrast, the standard deviation of the FT method remains positive even as we approach the free boundary, as explained above.

\begin{figure}[h]
	\begin{center}
		\begin{minipage}[t]{0.45\textwidth}
			\includegraphics[width=\textwidth]{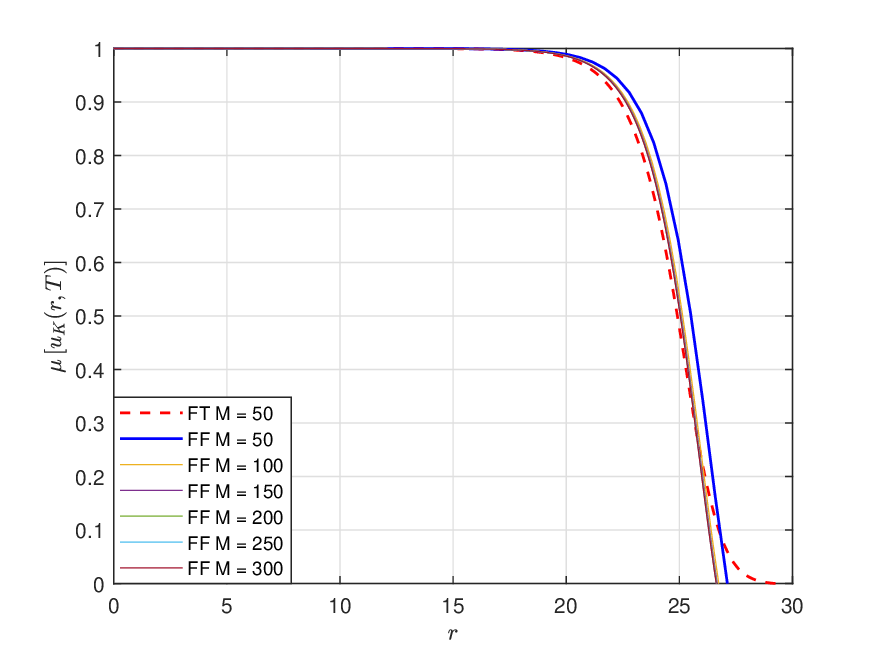}
		\end{minipage}
		\begin{minipage}[t]{0.45\textwidth}
			\includegraphics[width=\textwidth]{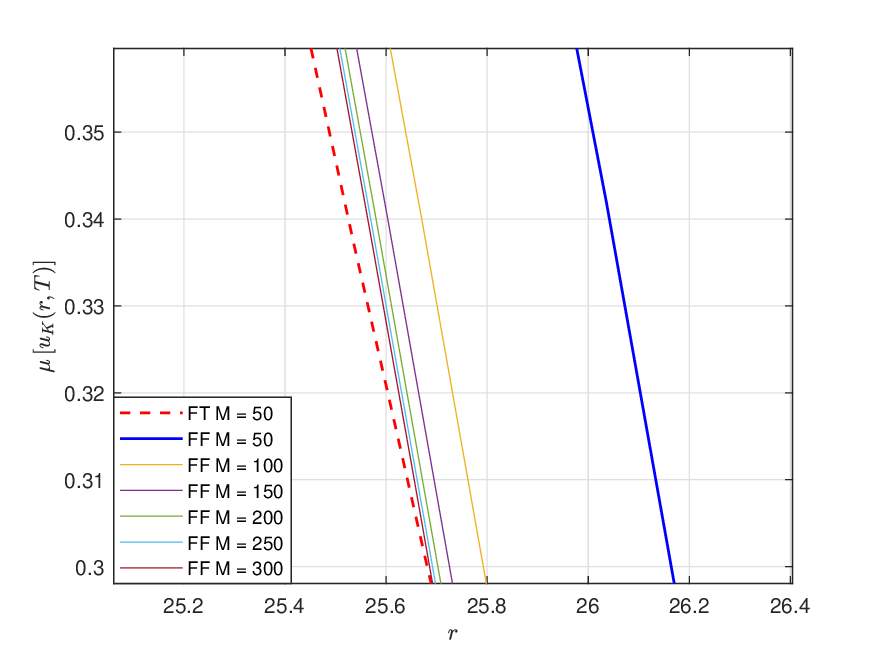}
		\end{minipage}
		\caption{Mean of the numerical solution s.p.'s for the  problem \eqref{eq:PDE}--\eqref{eq:ConditionsCI}, comparing the FF method with different numbers of grid points $M$, and the FT method using $M=50$ grid points.The plot on the right shows a zoomed-in view.}
		\label{fig:random_M}
	\end{center}
\end{figure}

Figure \ref{fig:random_M} displays the mean of the numerical solutions for problem \eqref{eq:PDE}--\eqref{eq:ConditionsCI}. It compares the FF method with various numbers of grid points $M$, against the FT method that uses $M=50$ grid points. Up to the inflection point, represented by the average of the free boundary as calculated by the FT method, the results exhibit convergence patterns that are similar to those observed in the previous deterministic case.

As demonstrated earlier, for short time intervals, both the FF and FT methods deliver similar results. Given its lower computational demands, we will utilize only the FF method in the next examples.


 We examine the spreading-vanishing dichotomy in the simplest scenario where both $\alpha$ and $\beta$ are constants. To facilitate this exploration, let's consider the following example.

\begin{example}\label{ex:spreading_vanishing}
Spreading-vanishing dichotomy.
\end{example}
In the logistic problem \eqref{eq:PDE}--\eqref{eq:BC} we set the default parameters in the constant case given in the first column of Table \ref{table:Comparison_parameters}.

\begin{figure}[h]
\begin{center}
\begin{minipage}[t]{0.45\textwidth}
	\includegraphics[width=\textwidth]{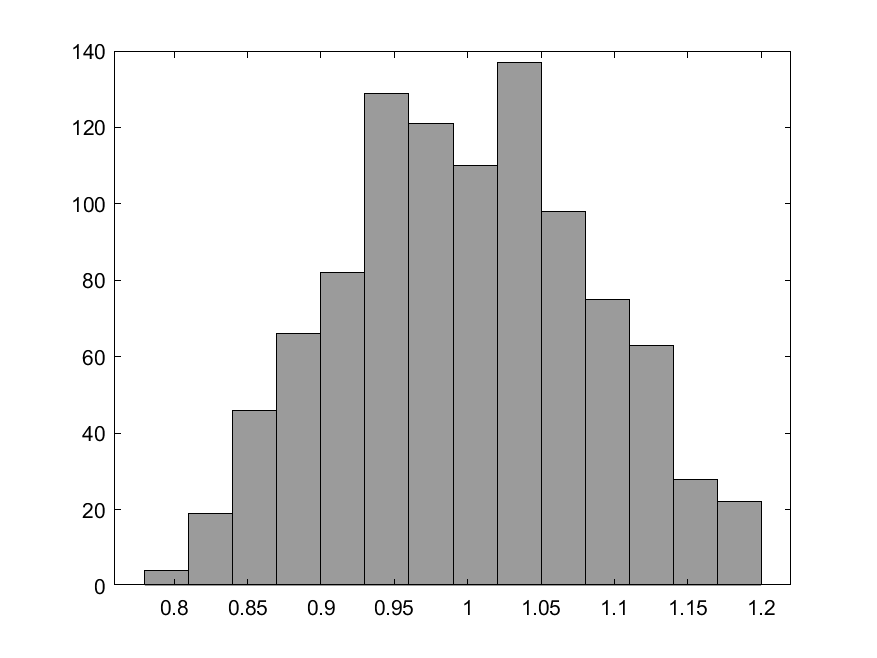}
\end{minipage}
\begin{minipage}[t]{0.45\textwidth}
	\includegraphics[width=\textwidth]{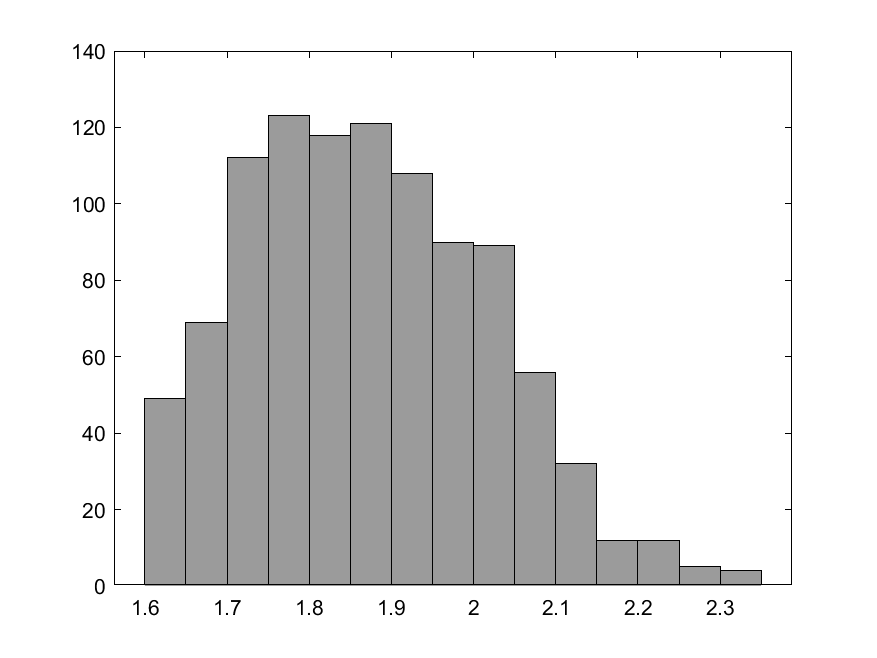}
\end{minipage}
\caption{Histograms of the sampled r.v.'s $D(\omega) \sim \mathcal{N}_{[0.8,\, 1.2]}(1, \, 0.1)$ (left) and $\mu(\omega) \sim \mathcal{B}e_{[1.6,\, 2.4]}(2, 4)$ (right).}
\label{fig:hist}
\end{center}
\end{figure}

For the Monte Carlo method, we set the number of samples as $K =10^3$, and for spatial nodes, we set $M = 50$. We also define $k =2\cdot 10^{-4}$ to ensure the stability of the FF method. The histograms of the sampled random variables $D(\omega)$ and $\mu(\omega)$, where $\omega$ is a member of the set $\Omega$, are illustrated in Figure \ref{fig:hist}.

\begin{figure}[h]
\begin{center}
    \includegraphics[width=0.45\linewidth]{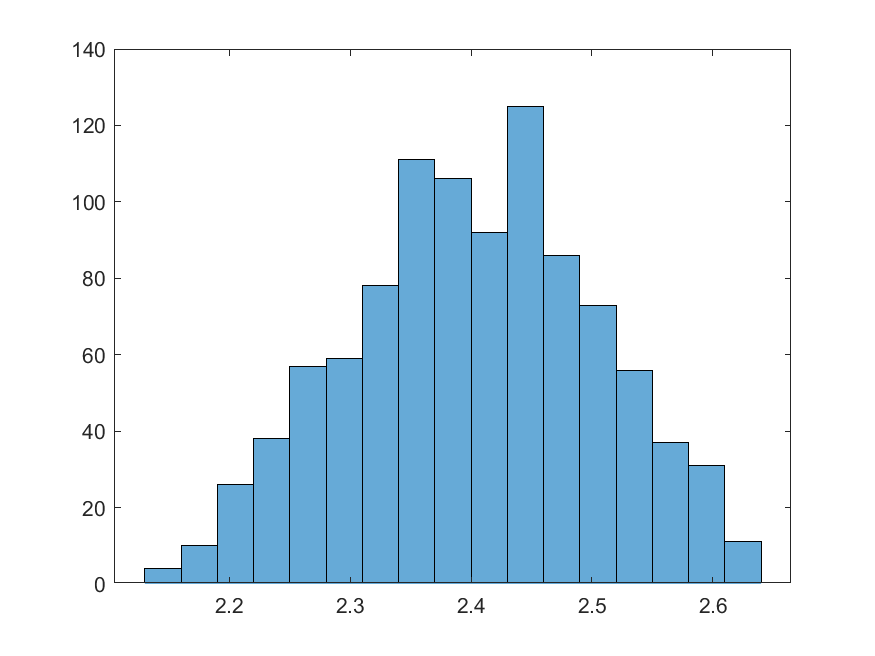}
\caption{Histogram of $R^*$ for Example \ref{ex:spreading_vanishing}. \label{fig:R}}
\end{center}
\end{figure}

Given that $\alpha$ is a constant, the spreading-vanishing boundary $R^*(\omega)$ can be determined analytically \citep{Casaban_deterministic}:
   \begin{equation}\label{eq:R_const}
  R^*(\omega) = r_0 \sqrt{\frac{D(\omega)}{\alpha}},\quad \omega \in \Omega,
  \end{equation}
  where $r_0$ is the first positive root of the Bessel function of the first kind, $r_0 = 2.40483$. Figure \ref{fig:R} presents the histogram of the distribution of $R^*(\omega)$, it is easy to notice that the histogram of $R^*$ replicated the shape of the histogram of $D(\omega)$.

  Formula \eqref{eq:R_const} shows that increasing $D$ leads to increasing $R^*$ for a fixed $\alpha$, hence

  \begin{equation}\label{eq:R_max}
     R^*_{\max} =  \max_{\omega \in \Omega} R^*(\omega) =  r_0 \sqrt{\frac{D_{\max}}{\alpha}}, \quad D_{\max}= \max_{\omega \in \Omega} D(\omega).
  \end{equation}

From \eqref{eq:R_max} and Theorem 1 of \citep{Casaban_deterministic}, the population spreading is guaranteed if $H_0 >  R^*_{\max}$. Since truncated normal distribution is used for generation of r.v. $D(\omega)$, then $0.8 \leq D(\omega) \leq 1.2$, for all realizations $\omega \in \Omega$, and $D_{\max} = 1.2$. For $\alpha = 1$, one gets 
\begin{equation}
    R^*_{\max} = r_0 \sqrt{\frac{D_{\max}}{\alpha}} = 2.6344.
\end{equation}

\begin{figure}[h]
\begin{center}
\begin{minipage}[t]{0.45\textwidth}
	\includegraphics[width=\textwidth]{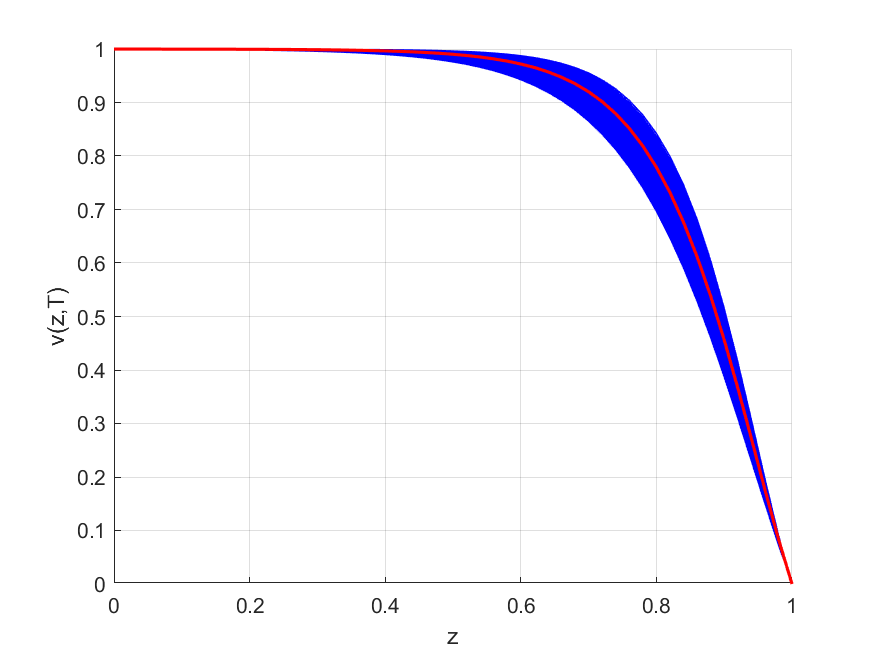}
\end{minipage}
\begin{minipage}[t]{0.45\textwidth}
	\includegraphics[width=\textwidth]{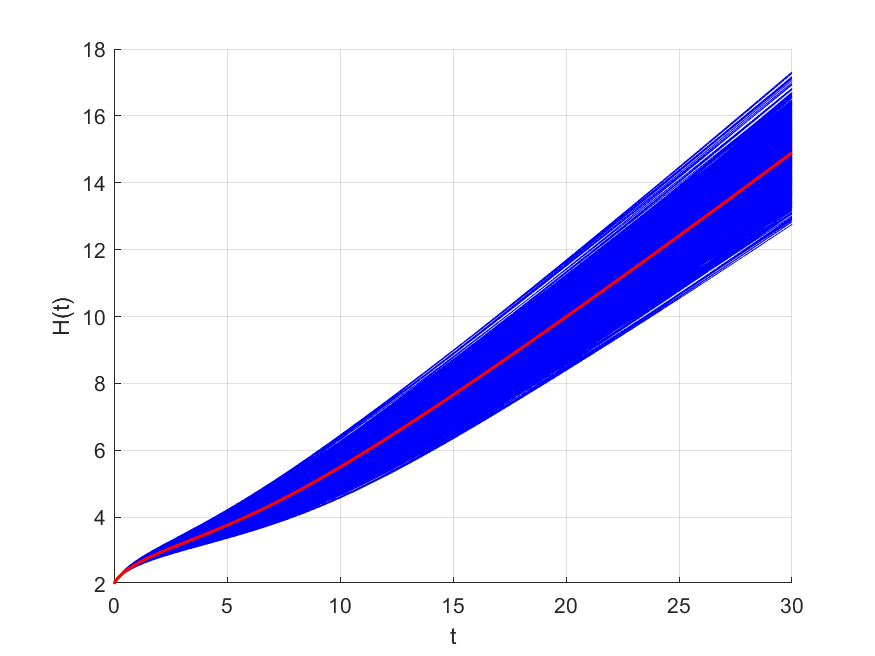}
\end{minipage}
\caption{Sampled numerical solution s.p.'s: transformed population distribution $v(z, t= T;\, \omega)$  (left) and the maximum population $v(z= 0, t;\, \omega)$, $\omega \in \Omega$, for $\mu(\omega) \sim \mathcal{B}e_{[1.6,\, 2.4]}(2, 4)$, the corresponding mean values are plotted in red.}
\label{fig:alpha_const_uH}
\end{center}
\end{figure}

Therefore, the population will spread for any choice of $H_0 > 2.6344$ for every realization, see Lemma \ref{lemma:R*}. In this specific example, $H_0 = 2$ is less than $R^*_{\max}$, which implies that spreading is not guaranteed and a further numerical examination of $\eta^*$ is necessary (see Theorem 1 in \citep{Casaban_deterministic} for more details). However, as illustrated in Figure \ref{fig:alpha_const_uH}, where the numerical solution and the free boundary s.p.'s are plotted along with their mean values, spreading occurs for all $\omega$ in the set $\Omega$.

\begin{figure}[h]
\begin{center}
\begin{minipage}[t]{0.45\textwidth}
	\includegraphics[width=\textwidth]{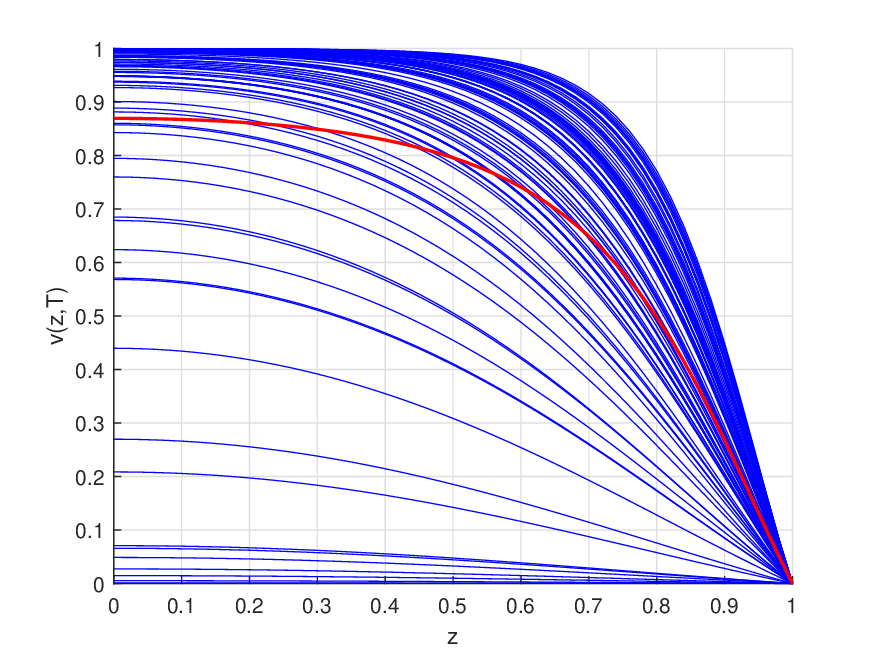}
\end{minipage}
\begin{minipage}[t]{0.45\textwidth}
	\includegraphics[width=\textwidth]{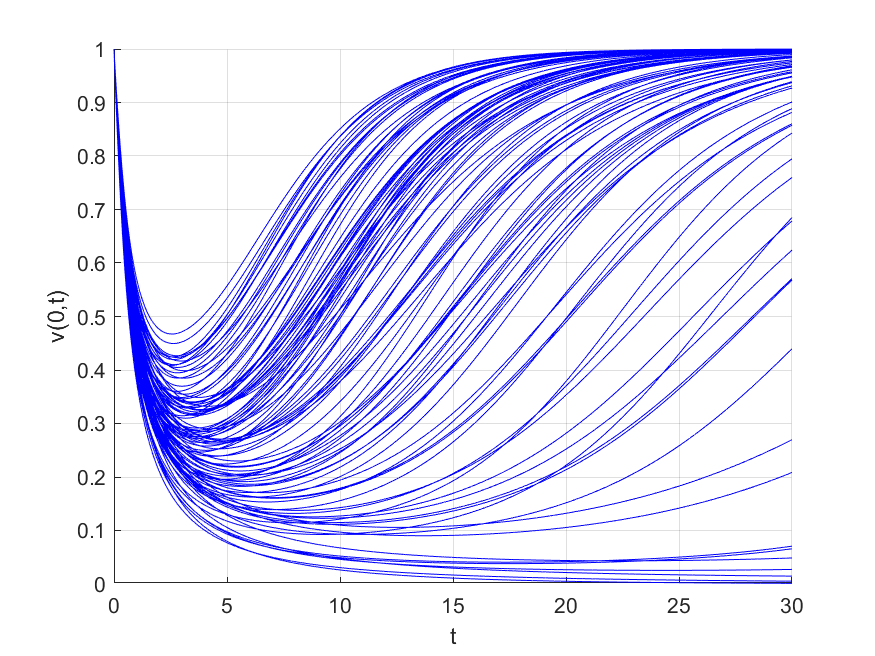}
\end{minipage}
\caption{Sampled numerical solution s.p.'s: transformed population distribution $V(z, t= T;\, \omega)$  (left) and the maximum population $V(z= 0, t;\, \omega)$, $\omega \in \Omega$, for $\mu(\omega) \sim \mathcal{B}e_{[0.2,\, 2.4]}(2, 4)$.}
\label{fig:alpha_const_SV_uH}
\end{center}
\end{figure}

Next, we set the range for $\eta(\omega)$ as $[0.2, 2.4]$ and reduce the number of samples to $K = 100$. This chosen range allows us to simulate scenarios where the potential for population extinction is possible given $H_0 = 2$. The sampled state variables for the population distribution and corresponding maximum population are displayed in Figure \ref{fig:alpha_const_SV_uH}. Observing the right plot, we note that only a few samples result in complete extinction, while the majority of samples indicate population spreading.

In the following example we shift our focus to the general case of $\alpha(r)$, $\beta(r)$. In this subsection, we numerically evaluate the convergence of the proposed random Finite Difference Method (RFDS-FF), along with its qualitative properties.

\begin{example} \label{ex:convergence}
	Numerical convergence.
\end{example}

\begin{figure}[h]
\begin{center}
\begin{minipage}[t]{0.45\textwidth}
	\includegraphics[width=\textwidth]{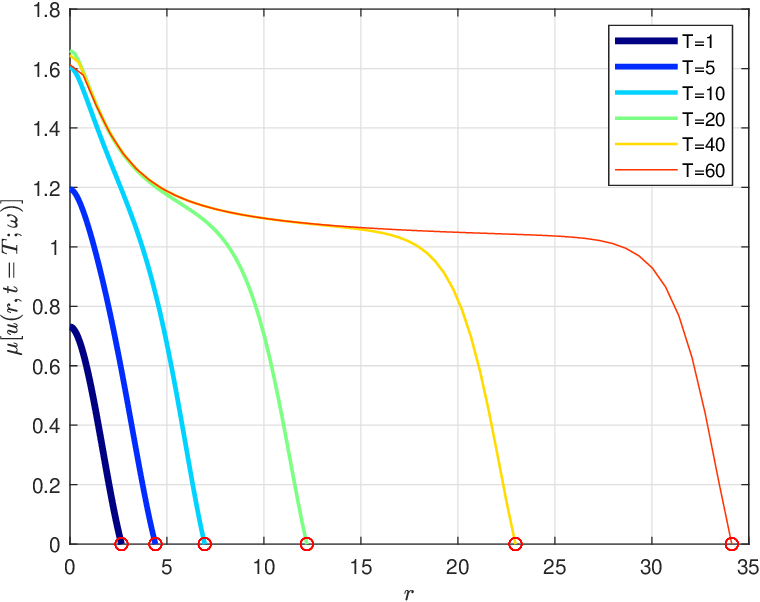}
\end{minipage}
\begin{minipage}[t]{0.45\textwidth}
	\includegraphics[width=\textwidth]{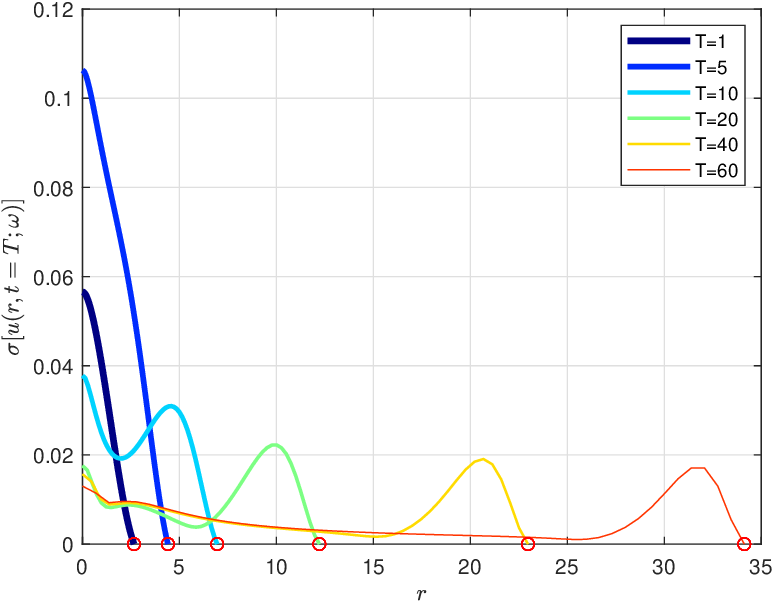}
\end{minipage}
\caption{Mean (left) and standard deviation (right) of the numerical solution s.p.'s for example \ref{ex:convergence}, for various time intervals $T$}
\label{fig:ex02_u_T}
\end{center}
\end{figure}

\begin{figure}[h]
\begin{center}
\begin{minipage}[t]{0.45\textwidth}
	\includegraphics[width=\textwidth]{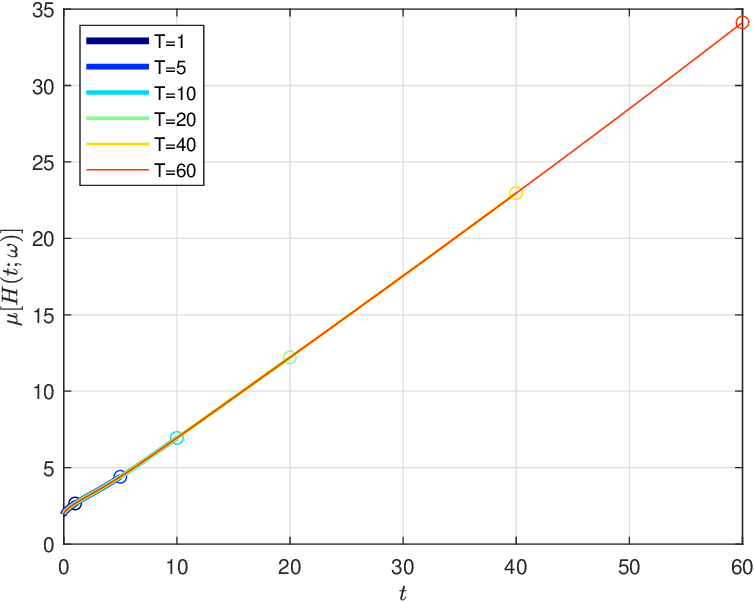}
\end{minipage}
\begin{minipage}[t]{0.45\textwidth}
	\includegraphics[width=\textwidth]{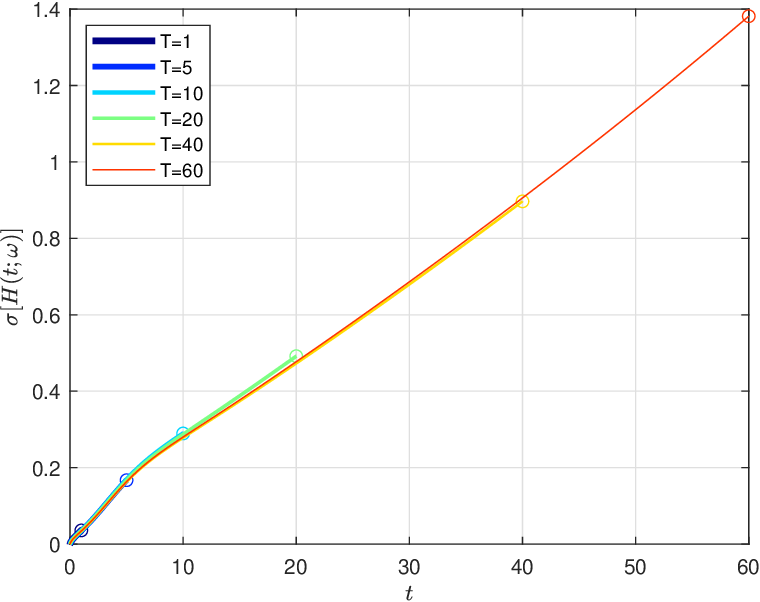}
\end{minipage}
\caption{Mean (left) and standard deviation (right) of the free boundary s.p.'s for example \ref{ex:convergence}, for various time intervals $T$}
\label{fig:ex02_H_T}
\end{center}
\end{figure}

Let us consider the random diffusive logistic model problem  \eqref{eq:PDE}--\eqref{eq:BC}  with the parameters from Table \ref{table:Comparison_parameters} (Variable case). 
The statistical moments of the numerical solution and free boundary s.p.'s for various $T$ values for the case of non constant $\alpha(r)$ and $\beta(r)$ are illustrated in Figures \ref{fig:ex02_u_T} and \ref{fig:ex02_H_T}, respectively. Note that in the left plot of Figure \ref{fig:ex02_u_T}, the red circles indicate the average free boundary at the corresponding time $t= T$ (this can be compared with the left plot of Figure \ref{fig:ex02_H_T}), demonstrating consistency between the two unknown s.p.'s.

\begin{figure}[h]
\begin{center}
\begin{minipage}[t]{0.45\textwidth}
	\includegraphics[width=\textwidth]{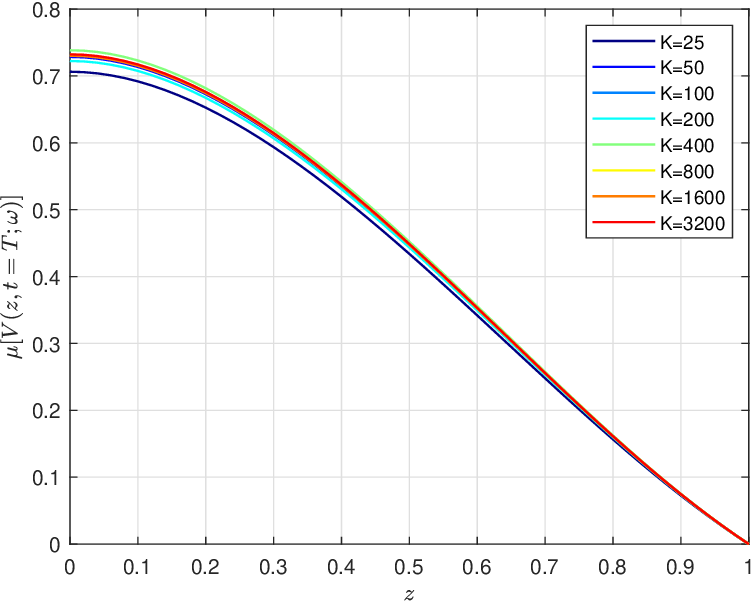}
\end{minipage}
\begin{minipage}[t]{0.45\textwidth}
	\includegraphics[width=\textwidth]{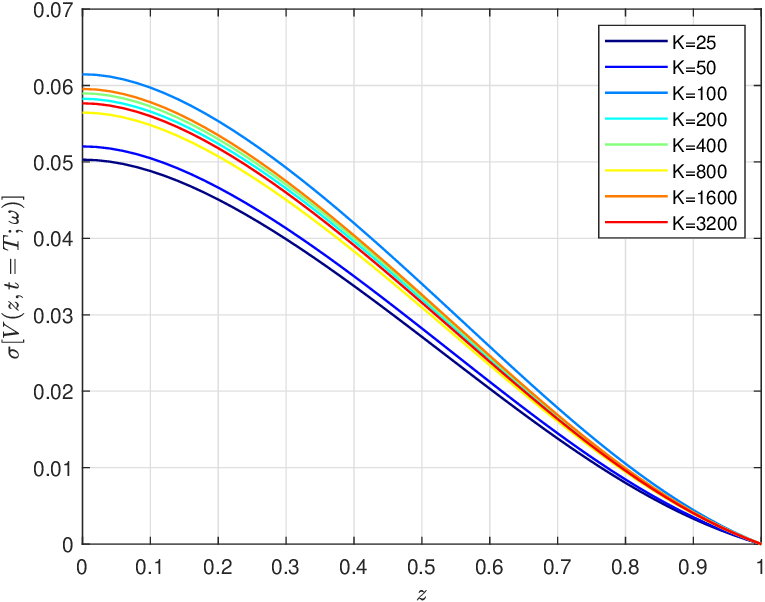}
\end{minipage}
\caption{Mean (left) and standard deviation (right) of the numerical solution s.p.'s for example \ref{ex:convergence}, for various number of samples $K$}
\label{fig:ex02_u}
\end{center}
\end{figure}

\begin{figure}[h]
\begin{center}
\begin{minipage}[t]{0.45\textwidth}
	\includegraphics[width=\textwidth]{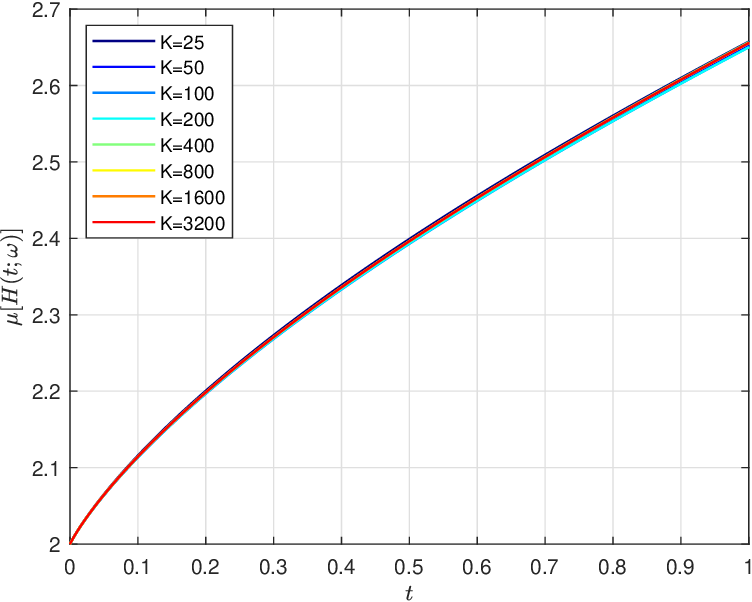}
\end{minipage}
\begin{minipage}[t]{0.45\textwidth}
	\includegraphics[width=\textwidth]{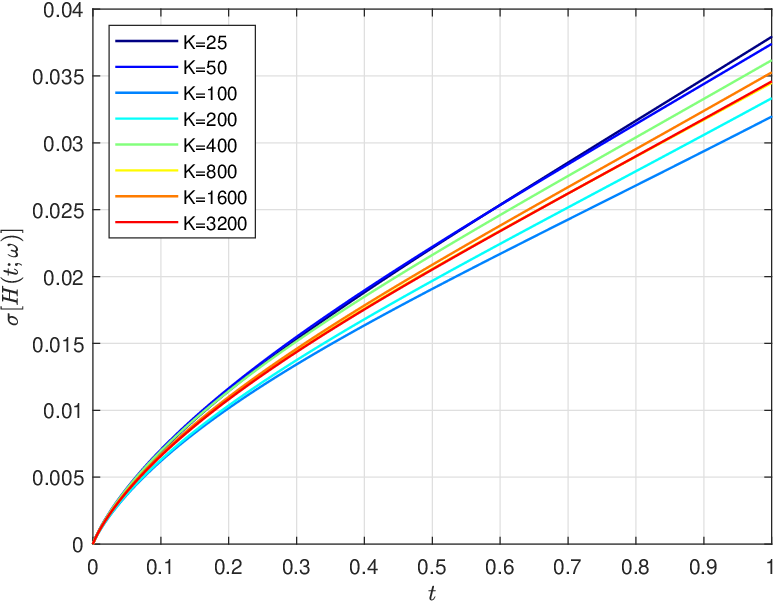}
\end{minipage}
\caption{Mean (left) and standard deviation (right) of the free boundary s.p.'s for example \ref{ex:convergence}, for various number of samples $K$.}
\label{fig:ex02_H}
\end{center}
\end{figure}

Firstly, we examine the numerical convergence of the Monte Carlo method. To minimize computational time, we set $T = 1$, which is sufficient for this convergence study. We establish the step-sizes by letting $M=50$, $N = 2000$ and we modify the number of samples for the Monte Carlo method by setting $K = 25\cdot 2^i$ for $0\leq i \leq 7$. The statistical moments for $V(z,t =T; \omega)$ and $H(t;\omega)$ are presented in Figures \ref{fig:ex02_u} and \ref{fig:ex02_H} respectively. In both instances, the mean values almost coincide for all considered $K$ values. The standard deviation of $V(z,T;\omega)$ is zero at $z = 1$ due to the fixed boundary condition, but it's positive at $z = 0$ because it's influenced by the random nature of the problem. The standard deviation of $H(t;, \omega)$ is zero at the initial moment, as $H_0$ remains constant for all $\omega \in \Omega$, but increases with time.

\begin{figure}[h]
\begin{center}
\begin{minipage}[t]{0.45\textwidth}
	\includegraphics[width=\textwidth]{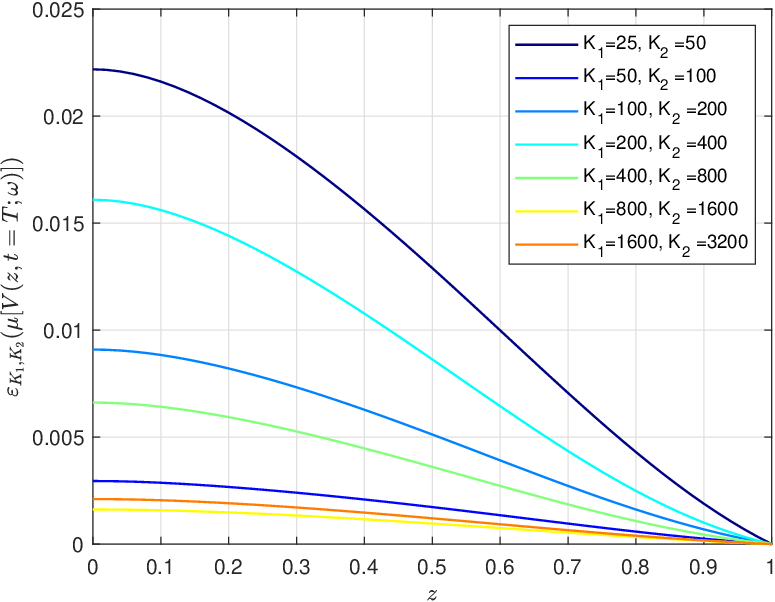}
\end{minipage}
\begin{minipage}[t]{0.45\textwidth}
	\includegraphics[width=\textwidth]{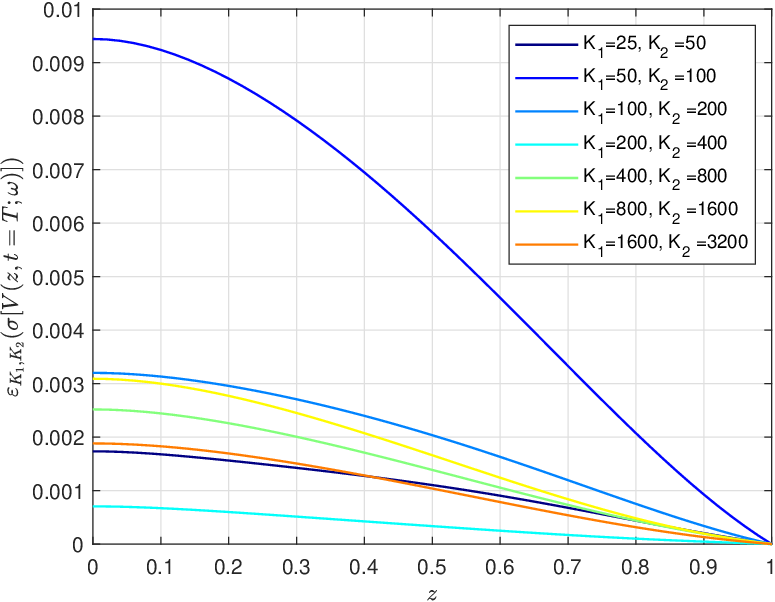}
\end{minipage}
\caption{Pairwise errors for the mean (left) and the standard deviation (right) of the numerical solution s.p.'s calculated by \eqref{eq:absdif_mu} and \eqref{eq:absdif:sigma}, respectively, for example \ref{ex:convergence}, for various pairs of $K$.}
\label{fig:ex02_absdif_u}
\end{center}
\end{figure}

\begin{figure}[h]
\begin{center}
\begin{minipage}[t]{0.45\textwidth}
	\includegraphics[width=\textwidth]{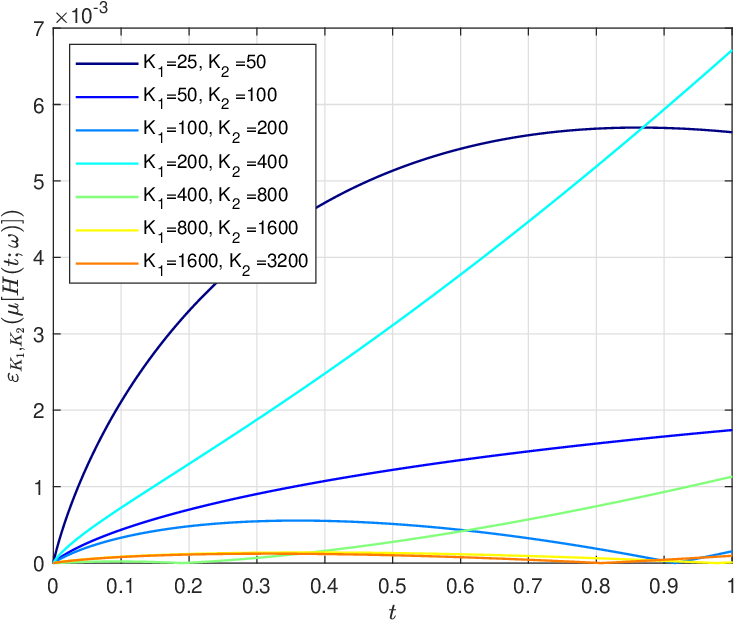}
\end{minipage}
\begin{minipage}[t]{0.45\textwidth}
	\includegraphics[width=\textwidth]{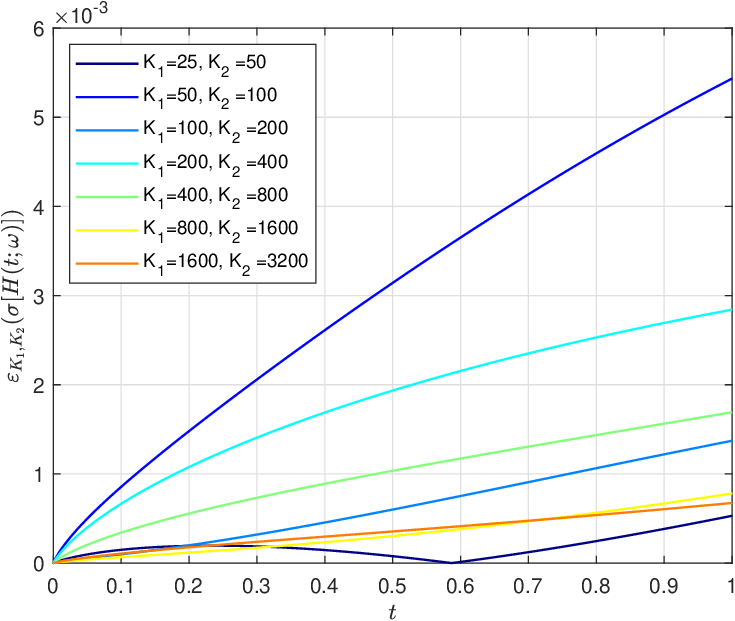}
\end{minipage}
\caption{Pairwise errors for the mean (left) and the standard deviation (right) of the free boundary s.p.'s calculated by \eqref{eq:absdif_mu} and \eqref{eq:absdif:sigma}, respectively, for example \ref{ex:convergence}, for various pairs of $K$.}
\label{fig:ex02_absdif_H}
\end{center}
\end{figure}

Visually from Figures \ref{fig:ex02_u} and \ref{fig:ex02_H}, the convergence of the Monte Carlo method is suspected but not obvious. Hence, in order to check the convergence, we introduce the following pairwise error measures
\begin{equation}\label{eq:absdif_mu}
    \varepsilon_{K_1, \, K_2}(\mu[u]) = \left| \mu_{K_1}[u] -\mu_{K_2}[u] \right|,
\end{equation}
\begin{equation}\label{eq:absdif:sigma}
    \varepsilon_{K_1, \, K_2}(\sigma[u]) = \left| \sigma_{K_1}[u] -\sigma_{K_2}[u] \right|,
\end{equation}
where $\mu_K[u]$ is the mean and $\sigma_K[u]$ is the standard deviation computed by the Algorithm \ref{algo:MC} for $K$ samples, $u$ stands for $V(z,T;\, \omega)$ or $H(t; \, \omega)$. These errors for various pairs of $K$ are plotted in Figures \ref{fig:ex02_absdif_u} and \ref{fig:ex02_absdif_H}.  The infinite norm of the errors is reported in Table \ref{table:ex02_absdif}. 

\begin{table}[h]
    \centering
    \begin{tabular}{|c|cc|cc|}
    \hline 
$ \left\{ K_1, \, K_2 \right\}$ &  ${\left\| \varepsilon_{K_1, \, K_2}(\mu[V]) \right\|}_{\infty}$ &          ${\left\| \varepsilon_{K_1, \, K_2}(\sigma[V]) \right\|}_{\infty}$  &  ${\left\| \varepsilon_{K_1, \, K_2}(\mu[H])\right\|}_{\infty}$ &          ${\left\| \varepsilon_{K_1, \, K_2}(\sigma[H]) \right\|}_{\infty}$ \\   \hline 
  $\left\{ 25 , \, 50\right\}$  & 2.2185e-02  & 1.7337e-03  & 5.6980e-03  & 5.3030e-04 \\ 
 $\left\{ 50 , \, 100 \right\}$  & 2.9472e-03  & 9.4415e-03  & 1.7405e-03  & 5.4339e-03 \\ 
 $\left\{ 100 , \, 200\right\}$  & 9.0891e-03  & 3.2002e-03  & 5.5588e-04  & 1.3724e-03 \\ 
 $\left\{ 200 , \, 400 \right\}$  & 1.6090e-02  & 7.0466e-04  & 6.7116e-03  & 2.8425e-03 \\ 
 $\left\{ 400 , \, 800 \right\}$  & 6.6110e-03  & 2.5169e-03  & 1.1311e-03  & 1.6915e-03 \\ 
 $\left\{ 800 , \, 1600 \right\}$  & 1.6162e-03  & 3.0901e-03  & 1.3896e-04  & 7.8092e-04 \\ 
 $\left\{1600 , \, 3200 \right\}$  & 2.1034e-03  & 1.8813e-03  & 1.2275e-04  & 6.7412e-04 \\ 
         \hline
         
    \end{tabular}
    \caption{Infinite norm of the pairwise error measure \eqref{eq:absdif_mu}--\eqref{eq:absdif:sigma} for the numerical solution and the free boundary s.p.'s. 
    \label{table:ex02_absdif}}
\end{table}

Now, let us study the numerical convergence of the FDM developed for the random FF method. For this purpose, we fix $K = 100$ and one of the step-sizes. By varying the second step size, and applying the analogous strategy as before, we calculate corresponding errors
\begin{equation}\label{eq:absdif_mu_h}
    \varepsilon_{M_1, \, M_2}(\mu[u]) = \left| \mu_{M_1}[u] -\mu_{M_2}[u] \right|,
\end{equation}
\begin{equation}\label{eq:absdif:sigma_h}
    \varepsilon_{M_1, \, M_2}(\sigma[u]) = \left| \sigma_{M_1}[u] -\sigma_{M_2}[u] \right|,
\end{equation}
\begin{equation}\label{eq:absdif_mu_k}
    \varepsilon_{N_1, \, N_2}(\mu[u]) = \left| \mu_{N_1}[u] -\mu_{N_2}[u] \right|,
\end{equation}
\begin{equation}\label{eq:absdif:sigma_k}
    \varepsilon_{N_1, \, N_2}(\sigma[u]) = \left| \sigma_{N_1}[u] -\sigma_{N_2}[u] \right|,
\end{equation}
where $M$ and $N$ are the number of spatial and temporal steps, respectively.

As in previous subsection, we set $T = 1$ that is enough for the convergence study. We fix $N = 10^6$, i.e., $k = 10^{-6}$ and vary $M$. The infinite norm of the pairwise error measure \eqref{eq:absdif_mu_h}--\eqref{eq:absdif:sigma_h} is reported in Table \ref{table:ex02_absdif_h}. For fixed $M = 50$ and various number of time-levels $N$, the infinite norm of the pairwise error measure \eqref{eq:absdif_mu_k}--\eqref{eq:absdif:sigma_k} is reported in Table \ref{table:ex02_absdif_k}. 

\begin{table}[h]
    \centering
    \begin{tabular}{|c|cc|cc|}
    \hline 
$ \left\lbrace M_1, \, M_2\right \rbrace$ &  ${\left\| \varepsilon_{M_1, \, M_2}(\mu[V]) \right\|}_{\infty}$ &          ${ \left\|\varepsilon_{M_1, \, M_2}(\sigma[V]) \right\|}_{\infty}$  &  ${\left\| \varepsilon_{M_1, \, M_2}(\mu[H]) \right\|}_{\infty}$ &          ${\left\| \varepsilon_{M_1, \, M_2}(\sigma[H]) \right\|}_{\infty}$ \\   \hline 
  $\left\lbrace 25 , \, 50 \right \rbrace$  & 7.4195e-04  & 4.0950e-03  & 6.0827e-04  & 1.4842e-03 \\ 
 $\left\lbrace 50 , \, 100 \right \rbrace$  & 3.1425e-03  & 3.2622e-03  & 4.5430e-03  & 1.5186e-03 \\ 
 $\left\lbrace 100 , \, 200 \right \rbrace$  & 6.0896e-03  & 5.1833e-03  & 1.3195e-02  & 4.8846e-03 \\ 
 $\left\lbrace 200 , \, 400 \right \rbrace$  & 4.5607e-03  & 1.6348e-03  & 1.3152e-02  & 6.0986e-04 \\ 
 $\left\lbrace 400 , \, 800 \right \rbrace$  & 2.8761e-03  & 6.7160e-04  & 2.2548e-03  & 4.5142e-04 \\ 
         \hline
         
    \end{tabular}
    \caption{Infinite norm of the pairwise error measure \eqref{eq:absdif_mu_h}--\eqref{eq:absdif:sigma_h} for the numerical solution and the free boundary s.p.'s. 
    \label{table:ex02_absdif_h}}
\end{table}

\begin{table}[h]
    \centering
    \begin{tabular}{|c|cc|cc|}
    \hline 
$ \left\lbrace N_1, \, N_2\right \rbrace$ &  ${\left\| \varepsilon_{N_1, \, N_2}(\mu[V]) \right\|}_{\infty}$ &          ${\left\|\varepsilon_{N_1, \, N_2}(\sigma[V]) \right\|}_{\infty}$  &  ${\left\| \varepsilon_{N_1, \, N_2}(\mu[H]) \right\|}_{\infty}$ &          ${\left\| \varepsilon_{N_1, \, N_2}(\sigma[H]) \right\|}_{\infty}$ \\   \hline 
  $\left\lbrace 2500 , \, 5000 \right \rbrace$  & 5.7087e-04  & 7.9248e-03  & 5.7833e-03  & 4.2857e-04 \\ 
 $\left\lbrace 5000 , \, 10000 \right \rbrace$  & 2.1073e-03  & 4.8277e-03  & 1.5807e-04  & 2.2447e-03 \\ 
 $\left\lbrace 10000 , \, 20000 \right \rbrace$  & 1.6419e-03  & 6.8516e-03  & 3.6567e-03  & 3.8902e-03 \\ 
 $\left\lbrace 20000 , \, 40000 \right \rbrace$  & 2.9447e-03  & 8.4592e-03  & 4.7567e-04  & 1.6902e-03 \\ 
 $\left\lbrace 40000 , \, 80000 \right \rbrace$  & 5.3923e-03  & 5.2467e-04  & 3.2057e-03  & 3.8129e-04 \\  
         \hline
         
    \end{tabular}
    \caption{Infinite norm of the pairwise error measure \eqref{eq:absdif_mu_k}--\eqref{eq:absdif:sigma_k} for the numerical solution and the free boundary s.p.'s. 
    \label{table:ex02_absdif_k}}
\end{table}

\section{Conclusions}  \label{sec:Conclusions}

In this paper,  a random free boundary diffusive logistic model with radial symmetry is introduced considering parameters with a finite degree of randomness and a random unknown function describing the propagation front. Two different approaches, namely front-fixing (FF) and front-tracking (FT) methods, have been employed to handle the free boundary. For numerical solutions,  the corresponding  explicit RFDS are used and numerical analysis is performed to establish stability and positivity conditions for both methods. 

Using RFDS in the m.s. sense involves storage problems due to the iteration for computing the expectation and the variance of the solution s.p. To avoid these drawbacks, the Monte Carlo technique is applied to compute the statistical moments of the approximating s. p. solution and the stochastic moving boundary solution. Numerical examples illustrate the convergence of the Monte Carlo method and highlight qualitative properties such as stability, convergence, and spreading-vanishing analysis.

One of the objectives of present study consists of comparing the proposed methods and identifying their advantages and drawbacks, as well as specify their areas of applicability. The random FF method requires a suitable prefixed number of grid points for all time iterations and sample realizations of the random variables. In
contrast, the random FT method employs a varying number of grid points over time and across sample realizations. Additionally, comparing the computational cost of both methods  it is found that the random FF method is more efficient in terms of time and memory usage.  However, it is important to note that the random FF method may sacrifice accuracy due to the fixed boundary inverse transformation, resulting in a larger step-size in the original variables. This limitation can be overcome by using a larger number of grid points. In contrast, the random FT method maintains a consistent step-size, eliminating this hidden drawback.

In summary, the choice between the random FF and FT methods depends on the specific problem at hand. The random FF method is suitable for smaller time simulations when a fast and accurate result is desired. On the other hand, the random FT method proves to be a better option for solutions with a larger time horizon. By understanding the strengths and limitations of each method, researchers can select the most appropriate approach for their specific application.

\section*{Financial Support}

This work has been partially supported by the Spanish Ministry of Economy and Competitiveness MINECO through the project PID2019-107685RB-I00  and by the Spanish State Research Agency (AEI) through the project PDC2022-133115-I00.

\section*{Conflict of Interests Statement}
The authors have no conflicts of interest to disclose.

\section*{Ethics Statement}
This research did not require ethical approval.

\bibliographystyle{elsarticle-num}
\bibliography{population_bib}

\begin{thebibliography}{10}
\expandafter\ifx\csname url\endcsname\relax
  \def\url#1{\texttt{#1}}\fi
\expandafter\ifx\csname urlprefix\endcsname\relax\def\urlprefix{URL }\fi
\expandafter\ifx\csname href\endcsname\relax
  \def\href#1#2{#2} \def\path#1{#1}\fi

\bibitem{Fisher_1937}
R.~A. Fisher, The wave of advance of advantageous genes, Annals of Eugenics 7
  (1937) 355--369.

\bibitem{Kolmogorov_1937}
A.~Kolmogorov, N.~Petrovsky, S.~Piscounov, \'{E}tude de l' \'{e}quations de la
  diffusion avec croissance de la quantit\'{e} de mati\`{e}re et son
  application a un probl\`{e}me biologique, Bull. Univ. Moskou 1 (1937) 1--25.

\bibitem{Brauer2012}
F.~Brauer, C.~Castillo-Chavez, Mathematical Models in Population Biology and
  Epidemiology, Springer, 2012.

\bibitem{Acevedo2012}
M.~Acevedo, M.~Marcano, R.~J. Fletcher-Jr., A diffusive logistic growth model
  to describe forest recovery, Ecological Modelling 244 (2012) 13--19.
\newblock \href {https://doi.org/10.1016/j.ecolmodel.2012.07.012}
  {\path{doi:10.1016/j.ecolmodel.2012.07.012}}.

\bibitem{Zhou2014}
P.~Zhou, D.~Xiao, The diffusive logistic model with a free boundary in
  heterogeneous environment, Journal of Differential Equations 256(6) (2014)
  1927--1954.
\newblock \href {https://doi.org/10.1016/j.jde.2013.12.008}
  {\path{doi:10.1016/j.jde.2013.12.008}}.

\bibitem{Pagnini_2014}
G.~Pagnini, A.~Mentrelli, Modelling wildland fire propagation by tracking
  random fronts, Natural Hazards and Earth System Sciences 14 (2014)
  2249--2263.

\bibitem{Durrett2010}
R.~Durrett, J.~Foo, K.~Leder, J.~Mayberry, F.~Michor, Evolutionary dynamics of
  tumor progression with random fitness values, Theoretical Population Biology
  78(1) (2010) 54--66.
\newblock \href {https://doi.org/10.1016/j.tpb.2010.05.001}
  {\path{doi:10.1016/j.tpb.2010.05.001}}.

\bibitem{Muntean2023}
S.~Nepal, Y.~Wondmagegne, A.~Muntean, Analysis of a fully discrete
  approximation to a moving-boundary problem describing rubber exposed to
  diffusants, Applied Mathematics and Computation 442 (2023) 127733.
\newblock \href {https://doi.org/https://doi.org/10.1016/j.amc.2022.127733}
  {\path{doi:https://doi.org/10.1016/j.amc.2022.127733}}.

\bibitem{Vynnycky2023}
M.~Vynnycky, On boundary immobilization for one-dimensional {S}tefan-type
  problems with a moving boundary having initially parabolic-logarithmic
  behaviour, Applied Mathematics and Computation 444 (2023) 127803.
\newblock \href {https://doi.org/https://doi.org/10.1016/j.amc.2022.127803}
  {\path{doi:https://doi.org/10.1016/j.amc.2022.127803}}.

\bibitem{Miklavcic2023}
M.~Miklav\v{c}i\v{c}, Analytic and numeric solutions of moving boundary
  problems, Journal of Computational and Applied Mathematics 431 (2023) 115270.
\newblock \href {https://doi.org/https://doi.org/10.1016/j.cam.2023.115270}
  {\path{doi:https://doi.org/10.1016/j.cam.2023.115270}}.

\bibitem{Friedman2015}
A.~Friedman, Free boundary problems in biology, Philosophical Transactions of
  the Royal Society A: Mathematical, Physical and Engineering Sciences
  373~(2050) (2015) 20140368.
\newblock \href {https://doi.org/10.1098/rsta.2014.0368}
  {\path{doi:10.1098/rsta.2014.0368}}.

\bibitem{DAcunto2021}
B.~D'Acunto, L.~Frunzo, V.~Luongo, M.~R. Mattei, A.~Tenore, Free boundary
  problem for the role of planktonic cells in biofilm formation and
  development, Zeitschrift f\"ur angewandte Mathematik und Physik 72~(149)
  (2021).
\newblock \href {https://doi.org/10.1007/s00033-021-01561-3}
  {\path{doi:10.1007/s00033-021-01561-3}}.

\bibitem{He2023}
T.~He, N.~Mitsume, F.~Yasui, N.~Morita, T.~Fukui, K.~Shibanuma, Strategy for
  accurately and efficiently modelling an internal traction-free boundary based
  on the s-version finite element method: Problem clarification and solutions
  verification, Computer Methods in Applied Mechanics and Engineering 404
  (2023) 115843.
\newblock \href {https://doi.org/10.1016/j.cma.2022.115843}
  {\path{doi:10.1016/j.cma.2022.115843}}.

\bibitem{Lu2020}
J.~Lu, B.~Hu, Bifurcation for a free boundary problem modeling the growth of
  multilayer tumors with ecm and mde interactions, Mathematical Methods in the
  Applied Sciences 43~(6) (2020) 3617--3636.
\newblock \href {https://doi.org/10.1002/mma.6142}
  {\path{doi:10.1002/mma.6142}}.

\bibitem{Chadam1993}
J.~M. Chadam, Emerging Applications in Free Boundary Problems, Chapman and
  Hall/CRC, 1993.

\bibitem{Vromans2018}
A.~J. Vromans, A.~Muntean, F.~van~de Ven, A mixture theory-based concrete
  corrosion model coupling chemical reactions, diffusion and mechanics, Pacific
  Journal of Mathematics for Industry 10~(5) (2018).
\newblock \href {https://doi.org/10.1186/s40736-018-0039-6}
  {\path{doi:10.1186/s40736-018-0039-6}}.

\bibitem{Du2010Spreading-VanishingBoundary}
Y.~Du, Z.~Lin, Spreading-vanishing dichotomy in the diffusive logistic model
  with a free boundary, SIAM Journal on Mathematical Analysis 42 (2010)
  377--405.
\newblock \href {https://doi.org/10.1137/090771089}
  {\path{doi:10.1137/090771089}}.

\bibitem{Du2022}
Y.~Du, Propagation and reaction-diffusion models with free boundaries, Bulletin
  of Mathematical Sciences 12~(01) (2022) 2230001.
\newblock \href {https://doi.org/10.1142/S1664360722300018}
  {\path{doi:10.1142/S1664360722300018}}.

\bibitem{TTSoong}
T.~Soong, Random Differential Equations in Science and Engineering, Academic
  Press: New-York, USA, 1973.

\bibitem{CasabanCompanyJodar2021}
M.~C. Casab\'{a}n, R.~Company, L.~J\'{o}dar, Reliable efficient difference
  methods for random heterogeneous diffusion reaction models with a finite
  degree of randomness, Mathematics 206 (2021) 351--366.
\newblock \href {https://doi.org/10.3390/math9030206}
  {\path{doi:10.3390/math9030206}}.

\bibitem{Du2011Spreading-vanishingII}
Y.~Du, Z.~Guo, Spreading-vanishing dichotomy in a diffusive logistic model with
  a free boundary, {II}, Journal of Differential Equations 250 (2011)
  4336--4366.
\newblock \href {https://doi.org/10.1016/j.jde.2011.02.011}
  {\path{doi:10.1016/j.jde.2011.02.011}}.

\bibitem{Du2020}
S.~Liu, Y.~Du, X.~Liu, Numerical studies of a class of reaction-diffusion
  equations with {S}tefan conditions, International Journal of Computer
  Mathematics 97~(5) (2020) 959--979.
\newblock \href {https://doi.org/10.1080/00207160.2019.1599868}
  {\path{doi:10.1080/00207160.2019.1599868}}.

\bibitem{CasabanCompanyJodarMATCOM2023}
M.~C. Casab\'{a}n, R.~Company, L.~J\'{o}dar, Numerical difference solution of
  moving boundary random {S}tefan problems, Mathematics and Computers in
  Simulations 205 (2023) 878--901.
\newblock \href {https://doi.org/10.1016/j.matcom2022.10.026}
  {\path{doi:10.1016/j.matcom2022.10.026}}.

\bibitem{Casaban_deterministic}
M.~C. Casab\'{a}n, R.~Company, V.~N. Egorova, L.~J\'{o}dar, Qualitative
  numerical analysis of a free-boundary diffusive logistic model, Mathematics
  11 (2023).
\newblock \href {https://doi.org/10.3390/math11061296}
  {\path{doi:10.3390/math11061296}}.

\bibitem{VillafuerteBraumannCortesJodar2010}
L.~Villafuerte, C.~A. Braumann, J.-C. Cort\'{e}s, L.~J\'{o}dar, Random
  differential operational calculus: Theory and applications, Comput. Math.
  Appl. 59 (2010) 115--125.
\newblock \href {https://doi.org/10.1016/j.camwa.2009.08.061}
  {\path{doi:10.1016/j.camwa.2009.08.061}}.

\bibitem{Okendal}
B.~Oksendal, Stochastic Differential Equations. An Introduction with
  Applications, Springer-Verlag, 1998.

\bibitem{Azbelev}
N.~Azbelev, Impact of certain traditions on development of the theory of
  differential equations, Computers and Mathematics with applications 37 (1999)
  1--8.

\bibitem{Crank}
J.~Crank, Free and Moving Boundary problems, Oxford University Press, 1984.

\bibitem{Gupta}
R.~S. Gupta, D.~Kumar, A modified variable time step method for the
  one-dimensional stefan problem, Comp. Meth. Appl. Mech. Engng. 23 (1980)
  101--108.

\bibitem{Marshall_2}
G.~Marshall, C.~Rey, L.~Smith, M\'{e}todos de seguimiento de la interface para
  problemas unidimensionales de frontera m\'{o}vil {II}, Revista internacional
  de m\'{e}todos num\'{e}ricos para c\'{a}lulo y dise\~{n}o en ingenier\'{i}a 2
  (1986) 351--366.

\bibitem{Kroese2011}
D.~P. Kroese, T.~Taimre, Z.~I. Botev, Handbook of Monte Carlo Methods, John
  Wiley et Sons Ltd., 2011.

\end{thebibliography}

%
%
%
%
%
%
%

\end{document}